\documentstyle[12pt]{article} 
\newfont{\twlvmsb}{msbm10 scaled\magstep1}
\newfont{\ninemsb}{msbm9}
\newfont{\sixmsb}{msbm6}
\newfam\msbfam
\textfont\msbfam=\twlvmsb
\scriptfont\msbfam=\ninemsb
\scriptscriptfont\msbfam=\sixmsb
\catcode`\@=11
\def\Bbb{\ifmmode\let\next\Bbb@\else
  \def\next{\errmessage{Use \string\Bbb\space only in math mode}}\fi\next}
\def\Bbb@#1{{\Bbb@@{#1}}}
\def\Bbb@@#1{\fam\msbfam#1}
\newfont{\largeeufm}{eufm10 scaled\magstep4}
\newfont{\twlveufm}{eufm10 scaled\magstep1}
\newfont{\elveufm}{eufm10 at 11pt}
\newfont{\teneufm}{eufm10}
\newfont{\nineeufm}{eufm9}
\newfam\eufam
\textfont\eufam=\twlveufm
\scriptfont\eufam=\nineeufm
\def\frak{\ifmmode\let\next\frak@\else
\def\next{\errmessage{Use \string\frak\space only in math mode}}\fi\next}
\def\frak@#1{{\fam\eufam{{#1}}}}
\catcode`@=12
\def\1{\hbox{ 1\kern-.35em\hbox{1}}}
\newcommand{\Z}{{\Bbb Z}} 
\newcommand{\C}{{\Bbb C}} 
\newcommand{\Nr}{{\Bbb N}_r} 
\newcommand{\cP}{{\cal P}} 
\newcommand{\cT}{{\cal T}_q(\frak g)} 
\newcommand{\Tq}{{\cal T}_q(\frak g)} 
\newcommand{\Eq}{{\cal E}_q^{\frak k}} 
\newcommand{\Eg}{{\cal E}_q} 
\newcommand{\Hq}{{\cal E}_q^{\frak k}} 
\newcommand{\Cq}{{\cal C}_q({\frak g}_0)} 
\newcommand{\Ck}{{\cal C}_q({\frak k})} 
\newcommand{\Fq}{{\cal F}_q} 
\newcommand{\Oq}{{\cal O}_q} 
\newcommand{\g}{{\frak g}} 
\newcommand{\Uq}{U_q({\frak g})} 
\newcommand{\RUq}{U^{\Bbb R}_q({\frak g}_0)} 
 
\newcommand{\Ul}{U_q({\frak l})} 
\newcommand{\RUk}{U_q^{\Bbb R}({\frak k})} 
\newcommand{\Up}{U_q({\frak p})} 
\newcommand{\id }{\mbox{id}} 
\newcommand{\quea}{ quantized universal enveloping algebra  } 
\newcommand{\queas}{ quantized universal enveloping algebras  } 
\newcommand{\irrep}{ irreducible representation  } 
\newcommand{\irreps}{ irreducible representations  } 
\newcommand{\Mod}{ \mbox{\bf Mod}_q(\g) } 
\newcommand{\ba}{\begin{eqnarray}}
\newcommand{\na}{\end{eqnarray}}
\newcommand{\ban}{\begin{eqnarray*}}
\newcommand{\nan}{\end{eqnarray*}}
\newtheorem{lemma}{Lemma}
\newtheorem{proposition}{Proposition}
\newtheorem{theorem}{Theorem}

\newcommand{\ol}[1]{\overline{#1}}
\newcommand{\gc}{\mbox{$G^{\Bbb C}$}}
\newcommand{\lig}{\mbox{$\frak g$}}
\newcommand{\lih}{\mbox{$\frak h$}}
 
\newcommand{\lip}{\mbox{$\frak p$}}
\newcommand{\liq}{\mbox{$\frak q$}}
\newcommand{\lik}{\mbox{$\frak k$}}
\newcommand{\lil}{\mbox{$\frak l$}}
\newcommand{\liu}{\mbox{$\frak u$}}
\newcommand{\bliu}{\ol{\mbox{$\frak u$}}}
 
\newcommand{\cb}{\mbox{$\cal B$}}
\newcommand{\ce}{\mbox{$\cal E$}}
\newcommand{\co}{\mbox{$\cal O$}}
\newcommand{\cq}{\mbox{$\cal Q$}}
\newcommand{\ct}{\mbox{$\cal T$}}
\newcommand{\cv}{\mbox{$\cal V$}}
\newcommand{\cw}{\mbox{$\cal W$}}
\newcommand{\cl}{{\cal L}^2_q}
\newcommand{\br}{\Bbb R}
\newcommand{\nd}{\nabla}
\newcommand{\dbar}{\ol{\partial}}
\newcommand{\nn}[1]{(\ref{#1})}

\newenvironment{proof}[1]{\begin{trivlist} \item[] {\em #1\/}: }%
{\hfill $\Box$ \end{trivlist}}
\newcommand{\ca}{{\cal A}_q({\frak g})}

\newcommand{\Hom}{\mbox{Hom}_{\C}}
\newcommand{\op}{\mbox{\scriptsize op}}

\begin{document}

\title{\normalsize{\bf GEOMETRY Of QUANTUM HOMOGENEOUS VECTOR BUNDLES AND
REPRESENTATION THEORY OF QUANTUM GROUPS I}}
\author{A. R. Gover$^\dagger$ and R. B. Zhang$^\ddagger$\\ 
 \small $^\dagger$ Mathematical Sciences, Queensland University of Technology, 
Brisbane, Australia\\
\small $^\ddagger$ Department of Pure Mathematics, Unversity of Adelaide, 
Adelaide, Australia}
\date{}
\maketitle
\begin{abstract}
Quantum homogeneous vector bundles are introduced by a direct
description of their sections in the context of Woronowicz type
compact quantum groups. The bundles carry natural topologies
inherited from the quantum groups, and their sections furnish
projective modules over algebras of functions on quantum homogeneous
spaces. Further properties of the quantum homogeneous vector bundles
are investigated, and their applications to the representation
theory of quantum groups are explored. In particular, quantum
Frobenius reciprocity and a generalized Borel-Weil theorem
are established.
\end{abstract}
 
\vspace{1cm}

\section{\normalsize INTRODUCTION}
The seminal work of Manin \cite{Manin} and Woronowicz \cite{Woronowicz}
demonstrated that quantum groups play much the same role in noncommutative 
geometry as that played by Lie groups in classical geometry.
This fact has been intensively investigated by several schools of 
researchers in recent years, considerably advancing our understanding 
of the underlying geometry of quantum groups.  We refer to the  
articles \cite{Zumino} and \cite{Majid} for reviews of the current state 
of the area and for useful references to the subject.  

The present paper is the first of a series intending to develop 
a comprehensive theory of quantum homogeneous vector bundles 
determined by quantum groups of the Woronowicz type\cite{groups}, 
and to explore their applications 
in a geometrical representation theory of quantum groups. 
Various versions of quantum deformations of fibre bundles 
were proposed 
at the algebraic level (i.e., without any topology) in the 
literature \cite{Brzezinski, Durdevic}. 
We mention in particular reference \cite{Brzezinski} 
and subsequent research along a similar line, where the primary aim was to 
develop a version of deformed gauge theory.  Quantum homogeneous 
vector bundles, in comparison, have been less studied, although they are 
much more closely related to quantum groups, and have 
natural applications in representation theory.   

In fact,  quantum homogeneous vector bundles                    
provide the foundations for developing a geometrical representation
theory of quantum groups.
There has long been an important interplay
between geometry and representation theory in the context of     
classical Lie groups, e.g., 
the interaction between representation theory
and the Penrose transforms of twistor theory \cite{Beastwood}. 
We expect a similar interaction between representation
theory and geometry to carry over to the quantum case.

The main new results of the paper are contained in Sections 3 and 4
where we define quantum homogeneous vector bundles and study their
properties and applications. In particular theorem \ref{BW} is a key
result.  In subsection \ref{bundles}, we introduce quantum homogeneous
vector bundles by a direct description of their sections. These are
defined in terms of Woronowicz type compact quantum groups and their
associated quantum homogeneous spaces.  Our definition of quantum
homogeneous vector bundles is consistent with the general definition
of noncommutative vector bundles adopted in Connes' theory
\cite{Connes}.  There a noncommutative vector bundle is defined by its
space of sections, which is required to be a projective module of
finite type over the algebra of functions on the noncommutative
(virtual) base space.  These noncommutative geometrical structures
carry natural topologies inherited from those of the quantum groups.
The latter is discussed in subsection \ref{completions}.  Projectivity
of the space $\Hq(V)$ of sections of a quantum homogeneous bundle
induced from a $\Ul$-module $V$ is established in Theorem \ref{proj}
of subsection \ref{projectivity}.  In this subsection, it is also
shown that $\Hq(V)$ forms an induced module over $\Uq$, and a
co-module over the associated quantum group, in analogy with the
classical situation.  Several other classical results are shown to
admit quantum analogues. In particular a quantum version of Frobenius
reciprocity is established in section \ref{Frobenius}, while
proposition \ref{projective} asserts that if the inducing $\Ul$-module
is, in fact, the restriction of a $\Uq$-module, then the space of
sections $\Hq(V)$ is freely generated as a module over the algebra of
functions of the quantum homogeneous space, and this means that the
quantum homogeneous bundle determined by $\Hq(V)$ is trivial.
Finally, in section \ref{BWsec} a notion of `quantum holomorphic'
sections is established and an analogue of the Borel-Weil theorem is
established. The reader may note that the approach to the proof there
is easily adapted to yield a new proof of the classical Borel-Weil
theorem via representative functions and the Peter-Weyl theorem.

The organization of the remainder of the paper is as follows.  Section
2 introduces the notation and conventions while reviewing the basic
definitions of quantum groups and quantized universal enveloping
algebras.  It also summaries their main structural and representation
theoretical features.  While the material, for the most part, is not new, our
treatment of real forms, parabolic and reductive quantum subalgebras
of \queas, as well as integrals on quantum groups, should be of general
interest. The appendix provides a concise and elementary treatment of 
the classical theory of homogeneous bundles as relevant to the quantum
constructions and results as mentioned above. Results are established there
in a manner that should shed light on the corresponding arguments for
the quantum case and indicate the geometrical nature of our treatment
of the latter. 

\section{\normalsize {\bf QUANTUM GROUPS AND QUANTIZED UNIVERSAL 
                     ENVELOPING ALGEBRAS}} 
\subsection{Quantized universal enveloping algebras}\label{QUEA}
\subsubsection{Generalities}
Let $\g$ be a finite dimensional simple complex Lie algebra of rank
$r$.  We denote by $\Phi^+$ the set of its positive roots relative to
a base $\Pi=\{\alpha_i \ | \ i\in\Nr\}$, where $\Nr=\{ 1, 2, ...,
r\}$.  Define $E=\bigoplus_{i=1}^r{\Bbb R}\alpha_i$.  Let $(\ , \ ):
E\times E \rightarrow {\Bbb R}$ be the inner product 
induced by the Killing form of $\g$. Then the Cartan matrix 
$A$ of $\g$ is given by $A=\left( a_{i j}\right)_{i j=1}^r$, with 
$a_{i j}={ {2(\alpha_i, \alpha_j)}\over{(\alpha_i, \alpha_i)}}$.  
We will call $\lambda\in E$
integral if $\lambda_i={ {2(\alpha_i, \lambda)} \over{(\alpha_i,
\alpha_i)}}\in\Z$, $\forall i$, and integral dominant if
$\lambda_i\in\Z_+$, $\forall i$.  The set of integral elements of $E$
will be denoted by $\cP$, and that of the integral dominant elements
by $\cP_+$.

The Jimbo version  \cite{Jimbo} of the   \quea  $\Uq$ is defined to be the
unital associative algebra over $\Bbb C$,  generated by $\{ k^{\pm 1}_i, 
e_i, f_i\ | \ i\in\Nr\}$  subject to relations as below. Here 
{\small  $\left[ \begin{array}{r}
       s\\ t
        \end{array}\right]_q $ } 
\normalsize 
is the Gauss polynomial, and  $q_i=q^{(\alpha_i, \ \alpha_i)/2}$. The
relations are:
\ba 
k_i k_j = k_j k_i, \quad  k_i k_i^{-1} =1, \quad 
k_i e_j k_i^{-1} = q_i e_j, \nonumber\\
k_i f_j k_i^{-1} = q_i f_j, \quad  
{[}e_i, \ f_j]=\delta_{i j}{{ k^2_i - k_i^{-2}}
\over{q_i -q_i^{-1}}}, \nonumber\\
\sum_{t=0}^{1-a_{i j}} (-1)^t
\left[ \begin{array}{r r r}
        1&-&a_{i j} \\
           &t&
        \end{array}\right]_{q_i}
(e_i)^t e_j (e_i)^{1-a_{i j}-t}= 0 , \quad  i\ne j,  \nonumber\\
\sum_{t=0}^{1-a_{i j}} (-1)^t
\left[ \begin{array}{r r r}
        1&-&a_{i j} \\
           &t&
        \end{array}\right]_{q_i}
      (f_i)^t f_j (f_i)^{1-a_{i j}-t}= 0 , \quad i\ne j, \label{quea}
\na 
where $q$ in general is taken to be a complex parameter, which is
nonvanishing and not equal to 1.  However, in this paper {\em we will
assume that $q$ is real positive and different from $1$}.
This  restriction is required  in order for $\Uq$ to admit a Hopf 
$\ast$-algebra structure and for the Haar functional on 
the corresponding quantum group to be positive definite. 
 
As is well known, $\Uq$ has the structure of a Hopf  algebra.  
We take the following co-multiplication
\ban
\Delta(k_i^{\pm 1})&=&k_i^{\pm 1}\otimes k_i^{\pm 1},\\
\Delta(e_i)&=&e_i\otimes k_i + k_i^{-1}\otimes e_i, \\
\Delta(f_i)&=&f_i\otimes k_i + k_i^{-1}\otimes f_i.
\nan
The co-unit $\epsilon: \Uq \rightarrow {\Bbb C}$
 and antipode $S: \Uq \rightarrow \Uq$ are respectively given by
\ban
\epsilon(e_i) =\epsilon(f_i)=0,  
&\epsilon(k_i^{\pm 1})= \epsilon(1)=1,& \\ 
S(e_i)=- q_i e_i, & S(f_i)=-q_i^{-1}  f_i, 
& S(k_i^{\pm 1})=k_i^{\mp 1}.
\nan 
 
The representation theory of $\Uq$ is closely related to that of 
the corresponding simple Lie algebra $\g$, and we refer to the 
many books on the subject, e.g., \cite{Chari}, for details. 
Here we mention that 
{\em all finite dimensional representations are completely reducible}.
Thus the study of such representations reduces to analyzing the 
irreducible ones.  If $W_\omega (\lambda)$ is a finite dimensional 
irreducible 
left $\Uq$-module, then the action of the $k_i$ can be diagonalized.  
There also exists a unique vector $v_+$ (up to scalar multiples), 
called the highest weight vector of $W_\omega (\lambda)$,  such that 
\ban 
e_i v_+ =0, & \quad  k_i v_+ = \omega_i q^{(\lambda, \ \alpha_i)/2} v_+,\\
\omega_i \in\{1, -1\}, & \quad \lambda\in\cP_+,    
\nan 
and the module $W_\omega (\lambda)$ is uniquely determined by 
$\lambda$ and the $\omega_i$.  
The existence of the $\omega_i$ is a peculiarity of
the Jimbo form of \quea, which stems from the
following  algebra automorphisms  
\ban 
e_i\mapsto \sigma_i e_i, & \quad f_i \mapsto \sigma^\prime_i
f_i & \quad k_i \mapsto \sigma_i\sigma^\prime_i  k_i,\\  
 & \sigma_i,  \ \sigma^\prime_i\in\{1, -1\}.&
\nan

When $\omega_i=1$, $\forall i$, we denote $W_\omega (\lambda)$ 
by $W(\lambda)$.  In this case, a common eigenvector $w\in W(\lambda)$ 
of the $k_i$ necessarily satisfies 
\ban 
k_i w &=& q^{(\mu,\ \alpha_i)/2} w, 
\nan 
for some $\mu\in\cP$. We call $\mu$ the weight of $w$.  The 
maximum weight, relative to the simple root system $\Pi$, is
$\lambda$. This  
will be referred to as the highest weight of $W(\lambda)$.    
When $\lambda\in\cP_+$, we will denote the lowest weight of 
$W(\lambda)$ by $\bar\lambda$, and define
\ban 
\lambda^\dagger&=&-\bar\lambda.
\nan 
Then $\lambda^\dagger$ is integral dominant,
and the dual module of $W(\lambda)$ has highest weight $\lambda^\dagger$.

Let $\Mod$ be the set of  
finite dimensional $\Uq$-modules,  
which is obviously closed  
under direct sum and direct product (with 
respect to the co-multiplication $\Delta$) 
In fact $\Mod$ forms a tensor category. The following points should be
observed:\\
{\em i). There is a one to one correspondence between the objects  
  of $\Mod$ and finite dimensional representations of the enveloping 
  algebra $U(\g)$ of $\g$;}\\ 
{\em ii). $W(\lambda)$,  $\lambda\in \cP_+$,  
  has the same weight space decomposition as 
  that of the irreducible $U(\g)$-module with highest weight $\lambda$.} 

\subsubsection{Real forms and parabolic subalgebras}
The \quea  $\Uq$ admits a
variety of Hopf $\ast$-algebra structures, namely, there exist anti
-involutions $\ast$ satisfying the following relation 
\ba \ast S \ast S &=& id_{\Uq}. \na 
Given an $\ast$-operation, we set 
\ba \theta &=& \ast S,   \label{theta} \na 
and call $\theta$ a quantum Cartan involution.  Let us define
\ba \RUq &=&\{ x\in\Uq\ | \ \theta ( x ) = x \}.  \na
It can be readily shown that
$\RUq$ defines a real associative algebra, which may be regarded as a
`real form' of $\Uq$. However, the restriction of $\Delta$ does not
lead to a co-multiplication for $\RUq$, thus $\RUq$ does not possess
a natural Hopf algebra structure.

Explicit examples of $\ast$-operations are the following, each 
specified by a choice of $\sigma_i,  \sigma'_i\in\{1, -1\}$: 
\ban
e_i^\ast =\sigma_i f_i, 
& f_i^\ast=\sigma'_i e_i,
& k_i^\ast=k_i^{\sigma_i \sigma'_i}.  
\nan
In this paper, we will be interested only in the {\em compact real form}
of $\Uq$. Thus, henceforth, we will assume that
$\RUq$ is defined by using the $\ast$-operation with
\ba
\sigma_i=\sigma'_i=1, \quad  \forall i.
\na
An important property  of this $\ast$-operation is that
every finite dimensional $\Uq$-module $W$ is unitary \cite{Gould}
in the sense that there exists a nondegenerate positive definite
sesquilinear form $(\ ,\ ): W\times W \rightarrow \C$ satisfying
\ba
( x v, \ w) &=& ( v, \ x^* w ), \quad \forall v, w\in W, \ x\in\Uq.
\label{sesquilinear}
\na

Denote by $\Cq$ the real vector space spanned by 
\ban 
X_i &=& e_i - q_i f_i, \\ 
Y_i &=&\sqrt{-1} ( e_i + q_i f_i ), \\ 
Z_i &=& \sqrt{-1} { {k_i - k_i^{-1}}\over{q_i - q_i^{-1}} }, \\ 
S_i &=& k_i + k_i^{-1} - 2, \quad \quad i\in\Nr.  
\nan 
Then $\RUq$ is generated by
$\Cq\cup\{\1_{\Uq}\}$.  Note that $\Cq$ vanishes under the co
-unit.  A further property is that $\Delta (\Cq) \subset
\Cq\otimes_{\Bbb R} \Uq + \Uq \otimes_{\Bbb R} \Cq$. That is, we have
the following result.
\begin{lemma}
$\Cq$ is a two-sided co-ideal of $\Uq$. 
\end{lemma}

Any element $a$ of the complexification of $(\RUq)^*$ naturally 
gives rise to a $\C$-linear functional on $\Uq$ ( regarded as  
the complexification of $\RUq$ ), by requiring that 
\ban 
a(x + \sqrt{-1} y) &=& a(x) + \sqrt{-1} a(y), \quad \forall x, y\in \RUq. 
\nan 
Conversely, we can restrict any linear functional on $\Uq$ to 
one on $\RUq$.  This identifies $\C\otimes_{\Bbb R} (\RUq)^*$ 
with $(\Uq)^*$.  In a similar way, one can easily
establish that there exists a one-to-one 
correspondence between complex 
representations of  $\RUq$ and complex 
representations of  $\Uq$.

Let us now consider parabolic and related subalgebras. 
Take a subset $\bf\Theta$ of $\Nr$.
Introduce the following sets of elements of $\Uq$:
\ban
{\cal S}_l&=&\{ k_i^{\pm 1}, i\in \Nr; \
\ e_j, \ f_j, \ j \in {\bf\Theta}\};\\
{\cal S}_{p}&=& {\cal S}_l \cup \{ e_j,
      j \in \Nr\backslash {\bf\Theta}\}.
\nan
Clearly ${\cal S}_l$ and 
${\cal S}_{p}$  generate Hopf subalgebras of $\Uq$, 
which we respectively denote by $\Ul$ and $U_q({\frak p})$. 
We call $\Ul$ a reductive  quantum subalgebra, and $U_q({\frak p})$
a parabolic quantum  subalgebra of $\Uq$, since, in the
classical limit, these Hopf subalgebras respectively reduce to
the enveloping algebras of  a reductive Lie subalgebra $\frak l$
and a parabolic subalgebra $\frak p$ of $\g$. 
We may replace the elements $e_j$, $j\in \Nr\backslash {\bf\Theta}$, 
by $f_j$ in ${\cal S}_{p}$, and the resulting set generates 
another Hopf subalgebra, which is the image of $\Up$ under 
the quantum Cartan involution. It also deserves the name of a parabolic 
quantum subalgebra. Results presented in the remainder of the paper can 
also be formulated using such {\em opposite} parabolic Hopf subalgebras.
It is important to observe that $\Ul$ is the invariant
subalgebra of $\Up$ under the quantum Cartan involution $\theta$.
For later use, we also define 
\ban \RUk &=& \Ul\cap \RUq. \nan  
Then $\RUk$ is a real subalgebra of $\RUq$, and its complexification 
is $\Ul$.  $\RUk$ is generated by $\1_{\Uq}$ and the set  
\ban 
\{ X_i, \ Y_i \ | \ i \in {\bf\Theta}\}\cup\{ Z_i, \ S_i \ | \ i \in \Nr\}. 
\nan  
We will denote by ${\cal C}_q(\frak k)$ the linear span of the elements of
this set. Then it can be easily shown that 
${\cal C}_q(\frak k)$ is a two-sided co-ideal of $\Uq$.  

Let $V_\mu$ be a finite dimensional irreducible $\Ul$-module.  Then
$V_\mu$ is of highest weight type.  Let $\mu$ be the highest weight
and $\tilde\mu$ the lowest weight of $V_\mu$ respectively.  We can
extend $V_\mu$ in a unique fashion to a $\Up$-module, which is still
denoted by $V_\mu$, such that the elements of ${\cal S}_{p}\backslash
{\cal S}_l$ act by zero.  It is not difficult to see that all finite
dimensional irreducible $\Up$-modules are of this kind.
 
Consider a finite dimensional irreducible $\Uq$-module
$W(\lambda)$, with highest weight $\lambda$  and lowest weight
$\bar\lambda$.  $W(\lambda)$ can be restricted  in a natural way
to a $\Up$-module, which 
is always indecomposable, but not irreducible in general.
It can be readily shown  that
\ba
\dim_{\C} \mbox{Hom}_{\Up} ( W(\lambda), \ V_\mu) &=&
    \left\{ \begin{array}{l l}
           1, &  \bar\lambda={\tilde\mu}, \\
           0, &  \bar\lambda\ne {\tilde\mu}.
     \end{array}\right.  \label{Hom}
\na

\subsection{Quantum groups}
\subsubsection{Quantum groups}
Roughly speaking a quantum group is the dual Hopf algebra of a 
\quea.   However, since $\Uq$ is 
infinite dimensional, considerable 
care needs to be exercised in defining the quantum group. 
The complication stems from the following well known fact. If 
$A$ is an infinite dimensional algebra, then the dual vector 
space $A^*$ in general does not admit a co-algebra structure.
The  way to get around the problem is to consider the so-call 
finite dual $A^0$ $\subset A^*$, which is defined by requiring that 
for any $f\in A^0$, $Ker f$ contains a two-sided ideal $\cal I$ 
of $A$ which is of finite co-dimension, i.e., 
$\dim A/{\cal I}<\infty$. 

A standard result of Hopf algebra theory states that if 
$A$ is a Hopf algebra with multiplication $m$, unit $\1_A$, 
co-multiplication $\Delta$, co-unit $\epsilon$ and antipode $S$, 
then the finite dual $A^0$, when it is not $0$, 
is also a Hopf algebra with a structure dualizing that of $A$. 
For any $ a, b\in A^0$, $x, y\in A$,  
the multiplication $m_0$ and co-multiplication $\Delta_0$ are respectively
defined by 
\ba 
\langle m_0(a\otimes b),\   x \rangle 
&=& \langle a\otimes b,\  \Delta(x) \rangle; \\ 
\langle \Delta_0 ( a ), \ x\otimes y \rangle
&=& \langle a, \ m(x\otimes y) \rangle;
\na  
the unit $\1_{A^0}$ and co-unit $\epsilon_0$  by  
\ban 
\1_{A^0}=\epsilon, \quad 
\epsilon_0 ( a ) = \langle a, \ \1_A \rangle;
\nan 
and the antipode $S_0$ by 
\ba 
\langle S_0(a), \ x\rangle 
&=& \langle a, \ S(x) \rangle.
\na 

Toward defining the quantum group dual to $\Uq$, consider the irreducible 
objects $W(\lambda)$, $\lambda\in\cP_+$ of $ \Mod$.  
For each $W(\lambda)$ of dimension $d_\lambda$, we choose a basis 
$\{ w_i^{(\lambda)}\ | \ i=1, 2, ...,d_\lambda\}$, 
which is arbitrary at this stage.   Let  $t^{(\lambda)}
= \left( t_{i j}^{(\lambda)}\right)_{i, j=1}^{d_\lambda}$, 
with $t_{i j}^{(\lambda)}$ being elements of $(\Uq)^*$ defined by 
\ba 
\sum_j t_{j i}^{(\lambda)}(x)  w_j^{(\lambda)} &=&
x  w_i^{(\lambda)},    \quad \forall x\in\Uq.         
\na 
We will also denote by $t^{(\lambda)}$ the \irrep of $\Uq$ associated 
with the module $W(\lambda)$ relative to the given basis, and 
call the $t_{j i}^{(\lambda)}$ matrix elements of the \irrep 
$t^{(\lambda)}$.

We denote by $\Pi(\lambda)$ the set of the weights of $W(\lambda)$. 
For each $i$, we defined $p_i^{\pm}$ 
$=\prod_{\nu\in\Pi(\lambda)}\left(k_i^{\pm 1} 
-q^{\pm (\nu, \alpha_i)/2}\right)$. 
Then $t^{(\lambda)}(p_i^{\pm})=0$, for all $i\in\Nr$.  
Let $\beta(\lambda)=\lambda-\bar{\lambda}$, where $\bar{\lambda}$ 
is the lowest weight of $W(\lambda)$.  If $u\in \Uq$ has the property 
that $k_i u k_i^{-1} = q^{(\gamma,\ \alpha_i)/2} u$, $\forall i\in\Nr$, 
and $\gamma>\beta(\lambda)$, or $\gamma < -\beta(\lambda)$, then  
$t^{(\lambda)}(u)=0$.  Let ${\cal I}_\lambda$ be the two-sided ideal 
of $\Uq$ generated by all such $u$ together with the $p_i^{\pm}$, 
$i\in\Nr$. Then it follows from  
the quantum analog of PBW theorem that ${\cal I}_\lambda$ has finite 
co-dimension, because $W(\lambda)$ is finite dimensional.  
Hence for all $\lambda\in\cP_+$, $t_{i j}^{(\lambda)}$ $\in$ 
$(\Uq)^0$.

The irreducibility of $W(\lambda)$ together with the  
so-called Burnside theorem of matrix algebras implies that 
$t^{(\lambda)}(\Uq)$ coincides with the entire algebra of 
$d_\lambda\times d_\lambda$ matrices. Hence for each $\lambda$ 
the $t_{i j}^{(\lambda)}$ are linearly independent. 
By considering the left action (\ref{circ}) of the central algebra of 
$\Uq$ on them, one can also easily convince oneself that the entire set 
$\{t_{i j}^{(\lambda)} | $ $i, j=1, 2, ...,d_\lambda$, $\forall$ 
$\lambda\in\cP_+\}$, is also linearly independent.  
Denote by $T^{(\lambda)}$ the subspace of $(\Uq)^0$ defined by  
\ban 
T^{(\lambda)} &=& \oplus_{i, j=1}^{d_\lambda} \C t_{i j}^{(\lambda)}. 
\nan 
Clearly this is independent of the choice of basis for $W(\lambda)$.  
Let 
\ba
\cT &=& \bigoplus_{\lambda\in\cP_+} T^{(\lambda)}, \label{group}
\na
where the direct sum is defined algebraically. 
The $\cT$ is essentially a quantum group of the kind introduced in 
\cite{Faddeev}.  It has the following important property  
\begin{proposition}
$\cT$ is a Hopf subalgebra of $(\Uq)^0$. 
\end{proposition}

This result is of course well known, but several aspects of it are  
worth mentioning.  Note that the multiplication of $(\Uq)^0$ is 
defined by the co-multiplication of $\Uq$.  Since the tensor product of
any two finite dimensional representations of $\Uq$ can be 
decomposed into a direct sum of finite dimensional \irreps,  
$\cT$ is indeed closed under multiplication. 
Let $R^{(\lambda) (\mu)}$ be the 
$R$-matrix associated with the two finite dimensional 
\irreps $t^{(\lambda)}$ and $t^{(\mu)}$, then 
\ban R^{(\lambda) (\mu)}_{1 2} t_1^{(\lambda)}t_2^{(\mu)} 
     &=& t_2^{(\mu)} t_1^{(\lambda)} R_{1 2}^{(\lambda) (\mu)}. 
\nan    
The co-product of each $t_{i j}^{(\lambda)}$ is easy to describe
explicitly. We have 
\ban 
\Delta_0(t_{i j}^{(\lambda)})
&=&\sum_{k=1}^{d_\lambda} t_{i k}^{(\lambda)}
\otimes t_{k j}^{(\lambda)}.    
\nan 
It is also useful to have an explicit characterization of the antipode. 
Introduce for $W(\lambda^\dagger)$ $=$ $(W(\lambda))^\ast$ 
the basis $\{ {\bar w}_i^{(\lambda)}\ |$ $i=1, 2, ..., d_\lambda\}$, 
which is  dual to the basis chosen for $W(\lambda)$ in the sense that 
${\bar w}_i^{(\lambda)} (w_j^{(\lambda)})$ $=$ $\delta_{i j}$.
Express the action of $\Uq$ on $W(\lambda^\dagger)$ by 
\ban 
x {\bar w}_i^{(\lambda)} =\sum_{j=1}^{d_\lambda} 
{\tilde t}_{ j i}^{(\lambda^\dagger)}(x)  
{\bar w}_j^{(\lambda)}, \quad x\in\Uq,   
\nan
for some ${\tilde t}_{ j i}^{(\lambda^\dagger)}$ in $\cT$. 
The natural action on the dual module is given by 
$(x v^\ast)(w)$ $=$ $v^\ast ( S(x) w)$, 
for any $v^\ast\in W(\lambda^\dagger)$,  $w\in W(\lambda)$ and
$x\in\Uq$. 
It then immediately follows that 
\ban 
\langle {\tilde t}_{ j i }^{(\lambda^\dagger)}, \ x \rangle 
&=&\langle S_0 ( t_{i j}^{(\lambda)}) , \ x\rangle,  
\quad \forall x\in\Uq, \\   
\mbox{i.e.}\quad  S_0(t_{i j}^{(\lambda)})&=& 
{\tilde t}_{ j i}^{(\lambda^\dagger)}.
\nan 
{ }From here on we will omit the subscript $0$ from $\Delta_0$ and $S_0$. 

In addition $\cT$ admits a natural anti-involution $\ast$-operation
giving it the structure of a Hopf $\ast$-algebra. The 
$\ast$-operation is defined by 
\ban \langle *(a), \ x \rangle &=&
\overline{\langle a, \ \theta(x) \rangle }, \quad \forall a\in\cT, \ 
x\in\Uq, \nan 
where $\theta$ is the quantum Cartan involution on
$\Uq$ defined by (\ref{theta}).  Simple computations can show that
this indeed gives rise to a $\ast$-operation for $\cT$.
It  takes the simplest form in a unitary basis of $\cT$, which 
we will introduce now. 
We assume that the basis $\{ w_i^{(\lambda)} \}$ 
of $W(\lambda)$ is orthogonal under the sesquilinear form 
(\ref{sesquilinear}). That is    
$ (w_i^{(\lambda)},  \ w_j^{(\lambda)} )$ $=$ $\delta_{i j}$.  
Define $\{ {\bar w}_i^{(\lambda)}\ | \ i=1, 2, ..., d_\lambda\}$ 
by ${\bar w}_i^{(\lambda)}( v )$ $=$ $(w_i^{(\lambda)},\ v  )$,
$\forall v\in W(\lambda)$,  
which form a dual basis for $(W(\lambda))^*$.  Note that  
$ (x {\bar w}_i^{(\lambda)} )( w_j^{(\lambda)} )$ $=$ 
$(\theta(x) w_i^{(\lambda)}, \ w_j^{(\lambda)} )$, 
which is equivalent to
$\langle {\tilde t}^{(\lambda^\dagger)}_{j i}, \ x\rangle$ 
$=$  $\overline{\langle t^{(\lambda)}_{j i}, \ \theta(x)\rangle}.$  
Thus in this basis, we have
\ba
\ast(t^{(\lambda)}_{i j}) &=&{\tilde t}^{(\lambda^\dagger)}_{i j}. 
\na  

\subsubsection{Quantum Haar measure}
The discussions of the last two subsections imply in particular 
that $\Tq$ satisfies the conditions of a CQG algebra in the sense of 
\cite{Koorwinder}.  Therefore the general theory of 
quantum Haar functionals of \cite{Ruan} \cite{Koorwinder} 
is applicable to $\Tq$.  Here we briefly treat the matter for the 
special case of $\Tq$. Our treatment  should be of general interest.

Let us begin by defining integrals on a general
Hopf algebra $A$ \cite{Montgomery}.
Let $A^*$ be its dual, which has a natural
algebra structure introduced by dualizing the co-algebraic
structure of $A$, although, as indicated above, 
$A^*$ does not admit a co-multiplication in general.   
An element $\int^l\in A^*$ is called a
left integral on $A$ if
\ban
f\cdot \int^l &=& \langle f, \1_A\rangle\ \int^l, \ \ \ \forall f\in A^*.
\nan
Similarly,  $\int^r\in A^*$ is called a right integral on $A$ if
\ban
\int^r\cdot f &=& \langle f, \1_A\rangle \ \int^r, \ \ \ \forall f\in A^*.
\nan
A straightforward calculation shows that the defining properties of
the integrals are equivalent to the following requirements
\ba
(\id\otimes \int^l)\Delta(x)=\int^l x, \quad 
(\int^r\otimes \id)\Delta(x)=\int^r x,\quad  \forall x\in A.
\na
where ${\rm id}$ is the identity map $A\rightarrow A$.

A normalised Haar measure $\int\in A^*$ on $A$ is an integral on $A$ which is
both left and right, and sends $\1_A$ to $1$, i.e.,
\ba
&(i).&  (\int\otimes \id) \Delta(x) = (\id\otimes \int)\Delta(x)
       = \int x, \quad \forall x\in A,\nonumber\\
&(ii).&  \int \1_A  = 1. \label{integral}
\na
When $A$ is a  Hopf $*$-algebra, we call a Haar measure
positive definite if $\int(x^* x)\ge 0$, and equality holds only
when $x=0$.

Now we go back to $\cT$. It is an entirely straightforward matter to 
establish the following result.
\begin{theorem}
The element $\int\in (\cT)^*$ defined by
\ba
\int \1_{\cT}=1; \quad  \quad \int t^{(\lambda)}_{i j}=0, 
\quad  0\ne \lambda\in \cP_+, \nonumber
\label{Haar} 
\na 
gives rise to a Haar measure on $\cT$.
\end{theorem}

Denote by $2\rho$ the sum of the positive roots of $\frak g$.  
Let $K_{2\rho}$ be the product of powers of $k_i^{\pm 1}$'s such that 
$K_{2\rho} e_i K_{2\rho}^{-1} = q^{(2\rho, \ \alpha_i)} e_i$, 
$\forall i$. Then it can be easily shown  that $S^2(x) = 
K_{2\rho} x K_{2\rho}^{-1}$, $\forall x\in\Uq$. Denote by  
$D_q(\lambda)$ $:=$ $tr\{t^{(\lambda)}(K_{2\rho})\}$ 
the quantum dimension of the irreducible $\Uq$-module  $W(\lambda)$.   
The Haar measure $\int$ satisfies the following properties.
\begin{lemma}\label{tresult}
\ba  
\int t^{(\lambda)}_{i j} {\tilde t}^{(\mu^\dagger)}_{r s} &=& 
{{t^{(\lambda)}_{s j}(K_{2\rho})}\over{D_q(\lambda)} } 
     \delta_{i r} \delta_{\lambda  \mu}, \nonumber \\
\int {\tilde t}^{(\lambda^\dagger)}_{i j} t^{(\mu)}_{r s} &=&
{{{\tilde t}^{(\lambda^\dagger)}_{ i r}(K_{2\rho})}\over{D_q(\lambda)} }
     \delta_{j s} \delta_{\lambda  \mu}.
\na  
\end{lemma}
{\em Proof}: The proof makes essential use of the fact that 
$\int$ is a left and right integral.  Look at the first 
equation.  The $\lambda \ne \mu$ case is easy to prove: 
the integral vanishes because   
the tensor product $W(\lambda)\otimes W(\mu^\dagger)$ 
does not contain the trivial $\Uq$-module. 
When   $\lambda = \mu$,  we introduce the notations  
\ban 
\phi_{i r; s j} = \int t^{(\lambda)}_{i j}  
                    {\tilde t}^{(\lambda^\dagger)}_{r s};
&\quad \Phi[s, j] = \left( \phi_{i r; s j}\right)_{i, r=1}^{d_\lambda};
&\quad \Psi[i, r]= \left( \phi_{i r; s j}\right)_{s, j=1}^{d_\lambda}.
\nan 
It is clearly true that 
$tr\left(\Psi[i, r]\right)$ $=$ $\delta_{i r}$. 

Note that corresponding to each $x\in\Uq$, there exists an  
${\tilde x}\in (\cT)^*$ defined by  
${\tilde x}(a) = \langle a, \ x\rangle$, $\forall a\in\cT$. 
The left integral property of $\int$ leads to 
\ban 
\epsilon(x)\phi_{i r; s j} &=& ({\tilde x}.\int)t^{(\lambda)}_{i j} 
{\tilde t}^{(\lambda^\dagger)}_{r s}\\
&=& \sum_{(x)} \sum_{i', r'} t^{(\lambda)}_{i i'}(x_{(1)}) 
{\tilde t}^{(\lambda^\dagger)}_{r r'}(x_{(2)}) \phi_{i' r'; s j},\\  
\mbox{i.e.}\quad  \epsilon(x) \Phi[s, j] 
&=& \sum_{(x)} t^{(\lambda)}(x_{(1)}) \Phi[s, j] t^{(\lambda)}(S(x_{(2)})), 
\quad \forall x\in\Uq.    
\nan 
Schur's lemma forces $\Phi[s, j]$ to be proportional to the identity matrix, 
and we have 
\ban 
\Psi[i, r]&=&\delta_{i r} \psi, 
\nan 
for some $d_\lambda \times d_\lambda$ matrix $\psi$. 
The right integral property of $\int$ leads to 
\ban \epsilon(x)\phi_{i r; s j} &=& 
(\int . {\tilde x} ) t^{(\lambda)}_{i j} 
{\tilde t}^{(\lambda^\dagger)}_{r s}\\ 
&=& \sum_{(x)} \sum_{j', s'} \phi_{i r; s' j'} t^{(\lambda)}_{j' j}(x_{(1)}) 
{\tilde t}^{(\lambda^\dagger)}_{s' s}(x_{(2)}).
\nan 
For $x=S(y)$, the equation is equivalent to
\ban 
\epsilon(y) \psi &=& \sum_{(y)} t^{(\lambda)}(K_{2\rho}) 
t^{(\lambda)}(y_{(1)}) t^{(\lambda)}(K_{2\rho}^{-1}) \psi 
t^{(\lambda)}(S(y_{(2)})). 
\nan 
Again by using Schur's lemma we conclude that $\psi$ is proportional 
to $t^{(\lambda)}(K_{2\rho})$. Since its trace is $1$, we have  
\ban \psi &=& { {t^{(\lambda)}(K_{2\rho})}\over{D_q(\lambda)} }. \nan 
This completes the proof of the first equation of the lemma. 
The second equation can be shown in exactly the same way. 

Given any $a=\sum_{\lambda\in\cP_+}\sum_{i, j=1}^{d_\lambda} 
           c^{(\lambda)}_{i j} t^{(\lambda)}_{i j}$,   we denote 
by $C^{(\lambda)}$ the matrix with entries $c^{(\lambda)}_{i j}$. 
Using the lemma we can easily show that 
\ban 
\int a^* a &=& \sum_{\lambda\in\cP_+} 
     tr\{{\tilde t}^{(\lambda^\dagger)}(K_{2\rho}) 
     C^{(\lambda)} {C^{(\lambda)}}^\dagger\}/D_q(\lambda), 
\nan    
which is clearly nonnegative, and vanishes only when $a=0$. 
We state this as a lemma.  
\begin{lemma}
The quantum Haar measure of $\cT$ is positive definite. 
\end{lemma} 

Note that the quantum Haar measure gives rise to a positive definite
sesquilinear form $(\cdot ,\cdot)_h$ for $\Tq$ defined by
\ban
(a,b)_h&=&\int a^\ast b , \quad a, b\in\cT.
\nan

\section{\normalsize{\bf QUANTUM HOMOGENEOUS VECTOR BUNDLES AND 
INDUCED REPRESENTATIONS}} 
\subsection{ Completions of ${\cal T}_q$ }\label{completions}
Existing treatments, in the literature, of quantum homogeneous spaces 
 are largely worked at the 
algebraic level, that is, without the introduction of topology. See for
example \cite{Schneider, Lakshmibai, Dijkhuizen}. In such cases the  
algebra of functions, over the quantum homogeneous space, is defined 
to be a subset 
of $\Tq$ satisfying certain homogeneity properties with respect to 
a two-sided co-ideal of $\Uq$.  
This is comparable to working with polynomials in a classical analysis
 situation. Thus, while such studies are instructive, it is ultimately
 unsatisfactory to remain in this purely algebraic setting.
To remedy this, we need to complete $\Tq$ in some way. 

Let $||\cdots||_h$ be the norm on $\cT$ determined by
\ban
||a||_h &=& \sqrt{ (a, a)_h }, \quad a\in\cT.
\nan
This equips $\Tq$ with the structure of a pre-Hilbert space.
Let us denote by $\cl$ the Hilbert space completion of $\cT$ 
in this norm.
Denote by $\cb(\cl)$ the bounded linear operators on $\cl$. 
Then the left regular representation of $\cT$ can be extended to 
the completion $\cl$, yielding a 
$\ast$-representation $\pi:\cT \to \cb(\cl)$ in 
the bounded operators $\cb(\cl)$.  
To prove this claim, note that for any $c\in\cT$,  
\ban 
(c a,b)_h&=&(a,c^\ast b)_h , \quad  a, b\in \cl, 
\nan 
if $|(c a,b)_h|<\infty$. 
Also observe that relative to a unitary basis 
(as discussed above in section \ref{QUEA}),
we have $\sum_k(t^{(\lambda)}_{ki})^\ast t^{(\lambda)}_{ki}=1$. 
Thus, for all $a\in \cl$,
\ban 
||a||_h^2&=&\sum_k(t^{(\lambda)}_{ki})^\ast t^{(\lambda)}_{ki}a,a)_h\\
 &=& \sum_k(t^{(\lambda)}_{ki}a, t^{(\lambda)}_{ki}a)_h\\
&\geq& ||t_{j i}^{(\lambda)} a||_h^2.  
\nan 
Therefore, the $t_{j i}^{(\lambda)}$ and their finite 
linear combinations  act on $\cl$ by  bounded 
linear operators.  
Let $||\cdot||$ be the operator norm on $\cb(\cl)$.  
For any $a\in\cT$,  $||\pi(a)||<\infty$.   
Thus, since $||\cdot||$ is a $C^\ast$-norm on
$\cb(\cl)$, the pull back of this under $\pi$ gives a $C^\ast$-norm
$||\cdot ||_{\op}$ on $\cT$ defined by 
\ba 
|| a ||_{\op}&=& \mbox{sup}\{ || a f||_h; \ f\in\cl, \ ||f||_h=1\}.
\na  
Finally, it is an elementary
exercise to check that the completion in this norm extends $\cT$ to a
unital $C^\ast$-algebra $\ca$.

The $C^\ast$-algebra $\ca$ qualifies as a compact quantum group of 
the Woronowicz type \cite{groups}, 
with $\Tq$ a dense subalgebra possessing the 
structure of a Hopf $\ast$-algebra.  However, we should note that 
it is not possible to extend the co-unit and antipode of $\Tq$ to 
continuous maps from the entire $\ca$ to appropriate spaces. 
Furthermore, an extension of the co-multiplication will necessarily 
maps $\ca$ continuously to some completion of $\ca\otimes\ca$ 
instead of the algebraic tensor product itself.   
 
As we have already mentioned, 
$\Tq$ satisfies all the conditions of a CQG algebra 
in the sense of \cite{Koorwinder}. Therefore, we could have followed 
that reference and completed $\Tq$ with respect to the largest 
$C^\ast$-seminorm, and consequently we would have  arrived 
at a different $C^\ast$-algebra.  However, 
our $\ca$ appears to be quite adequate for the purpose 
of studying quantum homogeneous bundles and investigating the 
representation theory of quantum groups, which is the main concern 
of our investigations.  

Let us introduce two types of actions of $\Uq$ on $\cT$. 
The first action will be denoted by $\circ$, which 
corresponds to the right translation in the classical theory 
of Lie groups.  It is defined by 
\ba
x\circ f&=&\sum_{(f)} f_{(1)} \ \langle f_{(2)}, \ x\rangle, 
\quad x\in\Uq, \ f\in\cT,  \label{circ}
\na
Straightforward calculations show that
\ban y\circ ( x\circ f ) &=& ( y x )\circ f;\\
(x\circ f) ( y ) &=& f(y x), \\
(id_{\cT}\otimes x\circ) \Delta(f) &=& \Delta( x\circ f ).  \nan
The other action, which corresponds to the
left translation in the classical  Lie group theory,
will be denoted by $\cdot$. It is defined by  
\ba
x\cdot f &=&\sum_{(f)} \langle f_{(1)}, \ S^{-1}(x) \rangle f_{(2)}. 
\na
It can be easily shown that  
\ban
(x\cdot f)(y)&=& f( S^{-1}(x) y ),\\
x\cdot(y\cdot f)&=& (x y)\cdot f, \quad  x,\ y\in \Uq, \ f\in \cT.
\nan
Furthermore, the two actions commute in the following sense  
\ban 
x\circ(y\cdot f)&=&y\cdot(x\circ f), \quad \forall x, y\in\Uq, \ f\in\cT.
\nan 

These actions can only be extended to subspaces of $\ca$. 
Any  $f\in \ca $ is the $n\to \infty$ limit of
a Cauchy sequence $\{ f_n\}$ with respect to the operator norm 
$||\cdot||_{\op}$, where each $f_n\in\Tq$.  Let $x\in\Uq$. 
We define 
\ban 
x\circ f=\lim_{n\to \infty} x\circ f_n,&\quad&
\mbox{if}\ || x\circ f_{n+m} - x\circ f_n||_{\op}\to 0, \quad n\to  \infty, \\
x\cdot f=\lim_{n\to \infty} x\cdot f_n,&\quad&
\mbox{if}\ || x\cdot f_{n+m} - x\cdot f_n||_{\op}\to 0, \quad n\to  \infty. 
\nan
Set 
\ba
{\cal E}_q:= \{ a\in\ca |  x\cdot a,\ x\circ a\in \ca,\  
                \ |a(x)|<\infty, \ \forall x\in\Uq\}.  
\na
The ${\cal E}_q$ clearly forms a subalgebra of $\ca$. 
In fact, for all $a,\ b\in\Eg$, 
we have 
\ban 
x\circ(a b) &=& \sum_{(x)} \{ x_{(1)}\circ a \} \{ x_{(2)}\circ b \}, \\ 
x\cdot(a b) &=& \sum_{(x)} \{ x_{(1)}\cdot a \} \{ x_{(2)}\cdot b \},  
\quad \forall x\in\Uq.
\nan 
We may regard $\Eg$ as the 
quantum analog of the algebra of smooth functions over the group.

\subsection{Quantum homogeneous spaces and 
  quantum homogeneous vector bundles}\label{bundles} 
Let us now turn
to the study of quantum homogeneous spaces.  As we will see shortly,
the well known fact in classical complex geometry, that any complex
analytic function on a compact complex manifold is a constant, also
holds in the analogous quantum setting.  Therefore, in the first
instance, we must work in a category of functions that has a richer
family of sections. This family should contain enough information to
capture the underlying geometrical aspects of the compact quantum
homogeneous spaces. On the other hand we want the class of functions
(and `bundle sections') to be closed under operations which generalize
classical differentiation. It is natural then to look for the quantum
analogs of algebras of smooth functions.  As in the classical case
(see section \ref{setup} of the appendix and in particular proposition
\ref{autohol}) this is most easily achieved by working in the `real
setting'. Thus we consider the compact real form of $\Uq$, and
regard $\cT$ as a subset of the complexification of $(\RUq)^*$.

Let us introduce the following definition   
\ba 
\Eq&:=&\left\{ f\in\Eg \ | \ x\circ f =\epsilon(x) f, 
\quad \forall x\in\RUk \right\}. \label{space} 
\na
Note that we may replace $\RUk$ by $\Ul=\C\otimes_{\Bbb R}\RUk$ 
in the above equation without altering $\Eq$.  
To investigate properties of $\Eq$, we 
consider the action of $\Ck$ on it. 
Recall that $\Ck$ generates the real subalgebra $\RUk$ of
$\RUq$. Also,  it is a two-sided co-ideal of $\Uq$  
and satisfies $\epsilon(\Ck)=0$. 
For any $a, \ b\in\Eq$, and $x\in\Ck$, we have 
\ban 
x\circ ( a b )&=&
\sum_{(x)} \{x_{(1)}\circ a\} \{ x_{(2)}\circ b\} =0. 
\nan 
Therefore $a b\in\Eq$, that is, 
\begin{center}
{\em $\Eq$ is a subalgebra of $\Eg$}.\\
\end{center}
We will show below that this non-commutative algebra is infinite 
dimensional.  We may regard it as the quantum analog of the algebra of
smooth functions on the homogeneous space $G_\C/P$, and  will refer to it 
as the algebra of functions on a (virtual) quantum 
homogeneous space.   
It is worth pointing out that in \cite{Dijkhuizen},  
more general quantum homogeneous spaces were considered, the definition 
of which was similar to (\ref{space}), but with $\RUk$ replaced by 
a two-sided  co-ideal of $\Uq$,  
which vanished under the co-unit and was $\theta$ invariant.
However,  the definition presented here is more 
suitable for the purpose of developing the representation 
theory of quantum groups.

Let $V$ be a finite dimensional module over $\Ul$,  which we will
also regard as a $\RUk$-module by restriction. 
We extend the actions $\circ$ and $\cdot$ of $\Uq$ on $\Eg$ trivially 
to actions on ${\cal E}_q\otimes V$: for any 
$\zeta=\sum_r f_r\otimes v_r\in \Eg\otimes V$   
\ban 
x\circ\zeta = \sum_r x\circ f_r\otimes v_r,&   
x\cdot\zeta = \sum_r x\cdot f_r\otimes v_r, & \quad   x\in\Uq. 
\nan  
We now introduce another definition, which will be of considerable 
importance for the remainder of the paper:  
\ban 
\Hq(V)&:=& \left\{ \zeta\in \Eg\otimes V \ 
| \ x\circ\zeta = (id_{\ca} \otimes S(x) ) \zeta, \ \forall x\in\RUk\right\}.
\nan  
Note that every $\zeta\in\Hq(V)$ 
satisfies 
\ban 
x\circ\zeta &=& (id_{\ca} \otimes S(x) ) \zeta, \quad  \ \forall x\in\Ul. 
\nan 

Consider the subspace 
$$\Fq(V):=\{ \Tq\otimes V\}\cap \Hq(V)$$ 
of $\Hq(V)$. 
Since the finite dimensional representations of $\Ul$ are completely 
reducible, the study of its properties reduces to the case 
when $V$ is irreducible. Let $V_\mu$ be a
finite dimensional irreducible $\Ul$ -module with highest weight $\mu$
and lowest weight $\tilde\mu$. 
Any element $\zeta$ $\in$  $\Fq(V_\mu)$
can be expressed in  the form
\ban
\zeta &=& \sum_{\lambda\in\cP_+}
\sum_{i, j} S(t_{ j i}^{(\lambda)})\otimes v_{i j}^{(\lambda)},
\nan
for some $v_{ i j}^{(\lambda)}\in V_\mu$.
Fix an arbitrary $\lambda\in\cP_+$. For any nonvanishing
$w\in W(\lambda)$, the following linear map is clearly surjective: 
\ban
\mbox{Hom}_{\C}(W(\lambda), \ V_\mu)\otimes w &\rightarrow& V_\mu,\\
\phi\otimes w&\mapsto& \phi(w). 
\nan
Thus there exist $\phi^{(\lambda)}_i$
$\in$ $\mbox{Hom}_{\C}(W(\lambda), \ V_\mu)$ such that
$v_{ i j}^{(\lambda)} = \phi^{(\lambda)}_i ( w_j^{(\lambda)} )$,
where $\{w_i^{(\lambda)}\}$ is the basis of $W(\lambda)$, relative to
which the \irrep  $t^{(\lambda)}$ of $\Uq$ is defined.
Now we can rewrite $\zeta$ as
\ban
\zeta &=& \sum_{\lambda\in\cP_+}
\sum_{i, j} S(t_{j i}^{(\lambda)})\otimes \phi_i^{(\lambda)}(w_j^{(\lambda)}).
\nan
The defining property of $\Fq(V_\mu)$ states that
\ban
\ell\circ\zeta&=& ( id_{\cT}\otimes S(\ell)) \zeta, \quad \forall \ell\in\Ul .
\nan
Thus  we have
$$
\sum_{\lambda\in\cP_+} \sum_{i, j, k} S(t_{k i}^{(\lambda)})
\otimes t_{j k }^{(\lambda)}(S(\ell)) \phi_i^{(\lambda)}(w_j^{(\lambda)})
=\sum_{\lambda\in\cP_+} \sum_{i, j} S(t_{j i}^{(\lambda)})\otimes S(\ell)
\phi_i^{(\lambda)}(w_j^{(\lambda)}) .
$$
Recalling that the $t_{k i}^{(\lambda)}$ are linearly
independent. It follows easily that the $S(t_{k i}^{(\lambda)})$ also
form a linearly independent set. So the above  is equivalent to
$$
\sum_{j} t_{j k }^{(\lambda)}(\ell) \phi_i^{(\lambda)}(w_j^\lambda)
= \ell \phi_i^{(\lambda)}(w_j^\lambda), ~~~ \forall \ell\in\Ul .
$$
This equation is precisely the statement that the $\phi_i^{(\lambda)}$
be $\Ul$-module homomorphisms,
\ban
\phi_i^{(\lambda)}&\in&\mbox{Hom}_{\Ul}\left(W(\lambda), V_\mu\right)
\subset \mbox{Hom}_{\C}\left(W(\lambda), V_\mu\right),
\quad \forall i.
\nan
Thus finding sections in $\Fq(V_\mu)$ is equivalent to finding, 
for all $\lambda\in\cP_+$, the homomorphisms $\phi^{(\lambda)}$ $\in$ 
$\mbox{Hom}_{\Ul}\left(W(\lambda), V_\mu\right)$. 
Note that each such homomorphism $\phi^{(\lambda)}$
determines $d_\lambda$ linearly independent sections
$$
\zeta^{(\lambda)}_i =  \sum_{j} S(t_{j i}^{(\lambda)})\otimes 
\phi^{(\lambda)}(w_j^{(\lambda)}). 
$$

Toward  constructing such homomorphisms we consider a couple of useful
observations. 
Note that if $W_1\to V_1$ and $W_2\to V_2$ are each $\Ul$-homomorphism 
then these induce a $\Ul$-homomorphism on the tensor product in 
the obvious manner
$$
 W_1\otimes W_2\to V_1\otimes V_2 .
$$
Now let $W(\lambda_1)$ and $W(\lambda_2)$ be irreducible
$\Uq$-modules of respective highest weights $\lambda_1$ and $\lambda_2$.
Let $V_{\mu_1}$ and $V_{\mu_2}$ be irreducible $\Ul$-modules
of the highest weights indicated. Then by explicit construction of maximal
weights one easily establishes the following:
\begin{lemma} \label{prodhom}
Suppose there are non-trivial $\Ul$-homomorphisms
$W(\lambda_1)\to V_{\mu_1}$
and $W(\lambda_2)\to V_{\mu_2}$. Then there is an induced non-trivial
$\Ul$-homomorphism
$$
W(\lambda_1+\lambda_2)\to V_{\mu_1+\mu_2}.
$$
\end{lemma}
 
Let us consider the case $\mu=0$, then  $\Fq(V_{\mu=0})=\Tq\cap\Eq$. 
We will show that this has an
infinite dimensional vector space of sections. Of course
there is a homomorphism from the trivial representation of $\Uq$,
$W(0)={\Bbb C}$, onto $V_0={\Bbb C}$. This gives the constant sections of
$\Tq\cap\Eq$.  Let $\gamma$ be the highest root of $\frak g$. 
Recall that $\Ul$ is reductive and there
are $N=r-|{\bf\Theta}|$ independent central elements in $\Ul$. 
Thus there are this many linearly independent $\Ul$-homomorphisms
$W(\gamma)\to {\Bbb C}$. 
As mentioned above each of these corresponds
to $d = \mbox{dim} ({\frak g})$  linearly independent sections.  
So the representation
$W(\gamma)$ determines $ N d$ 
linearly independent sections.
Further linearly independent sections may be obtained using lemma
\ref{prodhom}. For example there are $(m|N)$ 
(partition of $m$ into $\leq N$ parts) linearly independent homomorphisms
$W(m\gamma)\to {\Bbb C}$.  
It is easily verified that the $d
(m|N)$ sections so obtained are precisely the sections obtained by
taking $m$-fold products of the $d$ sections arising from the
homomorphisms $W(\gamma)\to V_\mu$.
We have proved the following lemma.
\begin{lemma}\label{pllenty} 
The algebra $\Eq$ is infinite dimensional. 
\end{lemma}

Now let us consider the case $\mu\neq 0$. It is an elementary exercise to
verify that $V_\mu$ is $\Ul$-isomorphic to a $\Ul$-irreducible
part of $W(\lambda')$, where $\lambda'$ is the dominant weight in the Weyl
group orbit of $\mu$. Thus there is a non-trivial $\Ul$-homomorphism
$$
W(\lambda') \to V_\mu ,
$$
and this  determines at least  $d_{\lambda'}$ linearly independent sections in
$\Fq(V_\mu)$.

Further linearly independent sections are obtained by
left and right multiplying with sections of $\Fq$. As above, the results of
such products may alternatively be constructed explicitly using lemma
\ref{prodhom} which promises a family of homomorphisms
$$
W(\lambda'+m\gamma)\to V_\mu ~~~~~m\in {\Bbb N}_+.
$$
Although we have fallen short of a classification of the sections in
$\Hq(V_\mu)$ we have established that $\Eq(V_\mu)$ is infinite
dimensional. This immediately leads to the following result.
\begin{proposition}\label{plenty}
If the weight of any vector of $V$ is $\Uq$-integral, then  
$\Eq(V)$ is an infinite dimensional vector space.  
\end{proposition}

$\Hq(V)$  provides a good candidate for the space of sections of  
a quantum vector bundle over the quantum homogeneous space 
corresponding  to $\Eq$.  We will discuss this further in the next section. 
Here we establish the following results.
\begin{theorem}
$\Hq(V)$ furnishes  \\
i). a two-sided  $\Eq$ module 
     under the multiplication of $\ca$;\\
ii). a  left $\Uq$-module under $\cdot$; and \\
iii). $\Fq(V)$ forms a left $\ca$ co-module under the co-action 
     $\omega = (\Delta\otimes id_V)$.  
\end{theorem}
{\em Proof}: Consider arbitrary elements $a\in\Eq$, 
$x\in\Uq$,  $p\in\RUk$, and $\zeta=$ $\sum_r f_r\otimes v_r$ 
$\in\Hq(V)$. 
The left and right actions of $\Eq$ 
on $\Hq(V)$ are respectively defined by 
\ban 
a\zeta=\sum_r a f_r\otimes v_r, &
\zeta a=\sum_r f_r a\otimes v_r.
\nan  
Now 
\ban 
p\circ(a\zeta) &=& \sum_{(p)}\{ p_{(1)}\circ a \} \{p_{(2)}\circ \zeta\}\\
               &=& \sum_{(p)} \epsilon( p_{(1)} ) a\{p_{(2)}\circ \zeta\}\\
               &=& a\{p\circ \zeta\} = (id_{\ca} \otimes S(p)) a\zeta;\\
p\circ(\zeta a)&=& \sum_{(p)} \{p_{(1)}\circ \zeta\} \{ p_{(2)}\circ a \}\\
               &=&\{p\circ \zeta\}a= (id_{\ca} \otimes S(p))\zeta a.
\nan 
This completes the proof of part i). Part ii) follows from 
\ban
p\circ  ( x\cdot\zeta )&=& x\cdot(p\circ\zeta)\\
&=&(id_{\ca} \otimes S(p)) ( x\cdot\zeta ), 
\nan 
while part iii) is confirmed by 
\ban 
(id_{\ca}\otimes  p\circ ) \omega(\zeta)
&=& (id_{\ca}\otimes  p\circ ) ( \Delta\otimes id_V)\zeta\\
&=&( \Delta\otimes id_V )(p\circ\zeta)\\ 
&=&\omega(id_{\ca}\otimes S(p))\zeta.            
\nan

Note that the left $\Uq$ action of ii) and left $\ca$ co-action  $\omega$ 
on $\Fq(V)$ are closely related.  Define a permutation map 
\ba 
P_{1 2 3}: \ca\otimes\ca\otimes V &\rightarrow& \ca\otimes V\otimes\ca
   \nonumber \\
      f_1\otimes f_2\otimes v&\mapsto& f_2\otimes v \otimes  f_1. 
\label{permutation}  
\na 
Then ${\tilde\omega}=P_{1 2 3}\omega$ 
defines a right $\ca$ co-action on $\Fq(V)$ which is 
dual to the left $\Uq$ action.  
We call $\Hq(V)$ an induced $\Uq$ module, and  also call 
$\Fq(V)$ an induced $\ca$ co-module. 

\subsection{Projectivity}\label{projectivity}
In classical differential geometry, the space $\cal H$ of sections of
a vector bundle over a compact manifold $M$ furnishes a module over
the algebra ${\cal A}(M)$ of functions.  It then follows from Swann's
theorem that this module must be projective and is of finite type,
namely, there exists another ${\cal A}(M)$-module ${\cal H}'$ such
that ${\cal H}$ $\oplus$ ${\cal H}^\prime$ is a finitely generated
free module, ${\cal A}(M)\oplus ... \oplus {\cal A}(M)$. (See also
theorem \ref{swannlike} of the appendix.) Conversely,
any projective module of finite type over ${\cal A}(M)$ is isomorphic
to the sections of some vector bundle over $M$.  This result is taken
as the starting point for studying vector bundles in noncommutative
geometry: one defines a vector bundle over a noncommutative space in
terms of the space of sections which is required to be a finite type
project module over a noncommutative algebra which is taken to be
the algebra of functions on the virtual noncommutative space.

Let us assume that all the weights of $V$ are integral, 
i.e., belonging to $\cP$.  In this case,   
$\Hq(V)$ will be called the space of 
sections of a quantum vector bundle over the quantum homogeneous 
space associated with $\Eq$.  
To justify this terminology, we need to show that 
$\Hq(V)$ is a projective module over $\Eq$.  
Let us first prove the following proposition. 
\begin{proposition}\label{projective}  
Let $W$ be a finite dimensional left $\Uq$-module, which we
regard as a left $\RUk$-module by restriction. 
Then $\Hq(W)$ is isomorphic to 
$\Eq\otimes W$ either as a left or right $\Eq$-module. 
\end{proposition}
{\em Proof}: We first construct the right $\Eq$-module 
isomorphism. 
Being a  left $\Uq$-module, $W$ carries a natural 
right $\ca$  co-module structure with the co-module action 
$\delta: W \rightarrow W\otimes\cT\subset W\otimes\ca$ 
defined for any element $w\in W$ by 
\ba 
\delta(w)(x) &=& x w, \quad  \forall x\in\Uq. 
\label{comodule} 
\na
Define the map $\eta: \Eg\otimes W \longrightarrow \Eg\otimes W$ 
by the composition of the maps
\ban 
\Eg\otimes W \stackrel{id\otimes\delta}{\longrightarrow} 
\Eg\otimes W\otimes\cT \stackrel{P^{-1}_{1 2 3}}{\longrightarrow} 
\cT\otimes \Eg\otimes W \longrightarrow \Eg\otimes W, 
\nan   
where the last map is the multiplication of $\ca$, 
and $P_{1 2 3}$ is the permutation map defined by (\ref{permutation}). 
Then $\eta$ defines a right $\Eq$-module isomorphism, 
with the inverse map given by the composition 
\ban
\Eg\otimes W \stackrel{id\otimes\delta}{\longrightarrow}
\Eg\otimes W\otimes\cT 
\stackrel{(S\otimes id\otimes id)P^{-1}_{1 2 3}}{\longrightarrow}
\cT\otimes \Eg\otimes W \longrightarrow \Eg\otimes W, 
\nan
where the last map is again  the multiplication of $\ca$.
It is not difficult to show that 
\ban 
x\circ\eta(\zeta)&=&\sum_{(x)}( id_{\ca}\otimes x_{(1)} ) 
                     \eta(x_{(2)}\circ \zeta),\\
x\circ\eta^{-1}(\zeta)&=&\sum_{(x)} 
( id_{\ca}\otimes S(x_{(1)}) ) \eta^{-1}(x_{(2)}\circ\zeta), 
\quad \forall \zeta\in\Eg\otimes W, \ \  x\in\Uq.
\nan 
Consider $\zeta\in\Hq(W)$. We have 
\ban 
p\circ\eta(\zeta)&=&  
   \sum_{(p)} ( id_{\ca}\otimes p_{(1)} )\eta(p_{(2)}\circ\zeta)\\
&=&\sum_{(p)}  ( id_{\ca}\otimes p_{(1)} S(p_{(2)}) )\eta(\zeta)\\ 
&=&\epsilon(p)\eta(\zeta), \quad  \forall \ p\in\Ul.
\nan 
Hence $\eta(\Hq(W))\subset \Eq\otimes W$. 
Conversely, given any $\xi\in\Eq\otimes W$, we have 
\ban 
p\circ\eta^{-1}(\xi)&=& 
=\sum_{(p)} ( id_{\ca}\otimes S(p_{(1)}) )\eta^{-1}(p_{(2)}\circ \xi )\\
&=&\sum_{(p)}( id_{\ca}\otimes \epsilon(p_{(2)}) S(p_{(1)}) )\eta^{-1}(\xi)\\   
&=& ( id_{\ca}\otimes S(p) ) \eta^{-1}(\xi),  \quad  \forall \ p\in\Ul.   
\nan 
Thus $\eta^{-1}(\Eq\otimes W)\subset{\cal E}_q(W)$. 
Therefore the restriction of $\eta$ provides the desired right 
$\Eq$-module isomorphism. 

The left module isomorphism is given by the restriction of 
$\kappa: \Eg\otimes W \rightarrow \Eg\otimes W$ defined by the 
composition of the following  maps
\ban
\Eg\otimes W \stackrel{id\otimes\delta}{\longrightarrow}
\Eg\otimes W\otimes\cT 
\stackrel{id\otimes(S^2\otimes id)P}{\longrightarrow}
\Eg\otimes \cT\otimes W \longrightarrow \Eg\otimes W,
\nan
where 
\ba 
P: W\otimes\ca&\rightarrow &\ca\otimes W, \nonumber \\
   w\otimes f &\mapsto& f\otimes w. \label{P}  
\na 
The inverse map $\kappa^{-1}$ is given by 
\ban
\Eg\otimes W \stackrel{id\otimes\delta}{\longrightarrow}
\Eg\otimes W\otimes\cT 
\stackrel{id\otimes(S\otimes id)P}{\longrightarrow}
\Eg\otimes \cT\otimes W \longrightarrow \Eg\otimes W. 
\nan

Let $V_\mu$ be a finite dimensional irreducible $\Ul$-module 
with highest weight $\mu$, which is integral with respect to $\frak g$.  
Then $V_\mu$ can always be embedded into an irreducible $\Uq$-module
$W(\sigma(\mu))$ with a $\frak g$ integral dominant highest weight  
$\sigma(\mu)$, where $\sigma$ is some element of 
the Weyl group ${\cal W}$ of $\frak g$. Such a $\sigma$ always exists, 
and belongs to the subgroup ${\cal W}^{\frak l}\subset {\cal W}$, 
which leaves invariant the set of the positive roots of $\frak l$. 
Since $\Ul$ is a reductive subalgebra of $\Uq$, all finite dimensional 
representations of $\Ul$ are completely reducible. Hence, $W(\sigma(\mu))$
can be decomposed into a direct sum of $\Ul$-modules:   
$W(\sigma(\mu))=V_\mu\oplus V_\mu^\bot$. 
Using the complete reducibility of finite dimensional $\Ul$-modules again, 
we conclude that
{if the weights of the finite dimensional $\Ul$ module $V$ 
are all integral with respect to $\Uq$, then there exist another 
$\Ul$ module $V^\bot$ and a  finite dimensional $\Uq$ module $W$ such that} 
\ban V\oplus V^\bot&=& W. \nan

It then immediately follows Proposition \ref{projective} that 
\ban 
\Hq(V) \oplus \Hq(V^\bot) 
&=& \Eq\otimes W, 
\nan
that is, 
\begin{theorem}\label{proj}
$\Hq(V)$ is projective and of finite type 
both as a left and right module over the algebra $\Eq$ of 
functions on the quantum homogeneous space. 
\end{theorem}

\subsection{Quantum Frobenius reciprocity}\label{Frobenius}
We have the following
quantum analog of Frobenius reciprocity.
\begin{proposition} Let $W$ be a $\Uq$ module, 
the restriction of which furnishes a $\RUk$ module in a natural way.
Then there exists a canonical isomorphism
\ba
\mbox{Hom}_{\Uq} ( W, \ \Hq(V)  )
&\cong& \mbox{Hom}_{\RUk} ( W, \  V),  
\na 
where $\Uq$ acts on the left module $\Hq(V)$ via the $\cdot$ action. 
\end{proposition}
{\em Proof}:  We prove the Proposition by explicitly constructing
the isomorphism, which we  claim to be  the linear map
\ban
F: \mbox{Hom}_{\Uq} ( W, \ \Hq( V ) ) &\rightarrow&
   \mbox{Hom}_{\RUk} ( W, \ V), \\
    \psi &\mapsto& \psi(1_{\Uq}),    
\nan
with the inverse map  
\ban
\bar{F}: \mbox{Hom}_{\RUk} ( W, \ V) &\rightarrow&
        \mbox{Hom}_{\Uq} ( W, \ \Hq( V ) ), \\
       \phi &\mapsto& \bar{\phi}=(S\otimes \phi) P\delta,   
\nan 
where $\delta: W\rightarrow W\otimes\Tq\subset W\otimes\ca$ is 
the right $\ca$ co-module action defined by (\ref{comodule}), 
and $P$ is the permutation map (\ref{P}).

As for $F$, we need to show that its image is contained in
$\mbox{Hom}_{\RUk} ( W, \ V)$.  Consider $\psi$ $\in$ 
$\mbox{Hom}_{\Uq} ( W, \ \Hq( V ) )$.  For any $p$ $\in$ $\RUk$ and 
$w$ $\in$ $W$,  we have 
\ban
p ( F \psi (w) ) &=&( (id_{\ca}\otimes p) \psi(w)) (\1_{\Uq})\\
                      &=& (S^{-1}(p)\circ\psi (w)) (\1_{\Uq}),  
\nan
where we have used the defining property of $\Hq( V )$.  
Note that 
\ban 
(S^{-1}(p)\circ\psi (w)) (\1_{\Uq})&=& (p\cdot\psi (w)) (\1_{\Uq}). 
\nan  
The $\Uq$-module structure of $\Hq( V )$ and the given condition 
that $\psi$ is a $\Uq$-module homomorphism immediately leads to 
\ban
p ( F \psi (w) ) &=& \psi(p w) (\1_{\Uq}) \\
                 &=& F \psi (p w), \quad\quad p\in \RUk, \  w\in W. 
\nan  

In order to show that ${\bar F}$ is the inverse of
$F$, we first need to demonstrate that the image $Im({\bar F})$
of ${\bar F}$ is contained in
$\mbox{Hom}_{\Uq} ( W, \ \Hq( V ) )$.
Note that $Im({\bar F})$ $\subset$  $\mbox{Hom}_{\C}
( W, \  \cT \otimes V )$. 
Some relatively simple manipulations lead to
\ban
(x\cdot{\bar\phi} (w) )  &=& {\bar\phi} ( x w ),\\
( p\circ {\bar\phi} (w) ) &=& (id_{\ca}\otimes S(p) ) 
                           {\bar\phi}(w), \ \ \ \
\forall x\in\Uq, \ \ p\in\RUk, \ \ w\in W.
\nan
Therefore,  $Im({\bar F})$ $\subset$
$\mbox{Hom}_{\Uq} ( W, \ \Hq( V ) )$.
Now we show that $F$ and ${\bar F}$ are inverse to each other.
For $\psi\in$ $\mbox{Hom}_{\Uq} ( W, \ \Hq( V ) )$, and
$\phi$ $\in$ $\mbox{Hom}_{\RUk} ( W, \ V)$, we have
\ban
(F \bar{F}  \phi)(w) &=& (\bar{F} \phi)(w) (1_{\Uq})\\
                             &=& \phi(w), \\
(\bar{F} F \psi)(w)(x)&=& ( F \psi) (S(x) w)\\
        &=& \psi (S(x) w) (1_{\Uq})\\
        &=&  (S(x)\cdot \psi (w)) (1_{\Uq})\\
        &=& \psi (w) (x), \ \ \ \ x\in\Uq, \ \ w\in W.
\nan
This completes the proof of the Proposition.

\section{\normalsize{\bf QUANTUM BOREL-WEIL THEOREM}}\label{BWsec}
Let $V_\mu$ be a finite dimensional irreducible $\Up$-module with
highest weight $\mu$ and lowest weight $\tilde\mu$. Recall that any two
norms on finite dimensional vector spaces determine the same topology.
Thus we may speak of convergence of a sequence in such
a space without reference to a particular norm.  Let us observe here 
that there is
a similar freedom for a certain class of norms on ${\cal
  E}_q\otimes_{\C} V_\mu$. To each basis $\{ v_r\}$ of $V_\mu$ we may
define a norm on ${\cal E}_q\otimes_{\C} V_\mu$  by 
\ban \zeta
=\sum_r f_r\otimes v_r, \quad\quad ||\zeta||^2=\sum_r ||f_r||_{\op}^2.
\nan 
It is easily verified that convergence in the norm corresponding
to one basis for $V_\mu$ implies convergence in all other norms
defined this way.  Thus, given $V_\mu$, we simply fix a basis and
define $||\cdot||$ to be the norm relative to that basis.

Recall the action of $\Uq$ on $\Eg\otimes V_\mu$. Since $V_\mu $ is a
$\Up$-module the following is a well defined subspace of 
${\cal E}_q\otimes_{\C} V_\mu$,
\ban \co_q(V_\mu) & :=& \left\{ \zeta\in {\cal E}_q
 \otimes V_\mu \ | \ p\circ \zeta =
(id_{\ca}\otimes S(p) )\zeta , \ \forall p\in\Up\right\}.
\nan
This may be regarded as the quantum analog of the space of 
holomorphic sections. 
Recall we use the notation $W(\lambda)$ to denote
the irreducible $\Uq$ module with highest weight $\lambda$.  We have
the following result.
\begin{theorem} \label{BW}
There exists the following $\Uq$ module isomorphism
\ba
\co_q(V_\mu) &\cong&\left\{\begin{array}{ l l }
                       W((-{\tilde \mu})^\dagger),& -{\tilde \mu}\in\cP_+,\\
                       0,& \mbox{otherwise}.
                       \end{array} \right.
\na
\end{theorem}
\begin{proof}{Proof of Theorem \ref{BW}} 
Let $\zeta$ $\in$ $\co_q(V_\mu)$.
Let $\{ \zeta_n\}$ be a  sequence in $\cT\otimes V_\mu$ such that 
$\zeta_n\to \zeta$ in the norm $||\cdot||$ described above.
Each $\zeta_n$ can be expressed in the form 
\ban 
\zeta_n &=& \sum_{\lambda\in\cP_+}
\sum_{i, j} S(t_{ j i}^{(\lambda)})\otimes v_{i j}^{(\lambda),n}, 
\nan
for some $v_{ i j}^{(\lambda),n}\in V_\mu$ ($i,j=1,\cdots d_\lambda$).
Arguing as in the proof of proposition \ref{plenty} (and lemma \ref{pllenty})
one concludes, for each $\lambda\in\cP_+$, that
there exist $\phi^{(\lambda),n}_i$ $\in$ $\mbox{Hom}_{\C}(W(\lambda), \ 
V_\mu)$ such that $v_{ i j}^{(\lambda),n} = \phi^{(\lambda),n}_i (
w_j^{(\lambda)} )$, where $\{w_i^{(\lambda)}\}$ is the basis of
$W(\lambda)$, relative to which the \irrep $t^{(\lambda)}$ of $\Uq$ is
defined.  Now we can rewrite $\zeta_n$ as 
\ban
\zeta_n &=& \sum_{\lambda\in\cP_+} \sum_{i, j} S(t_{j
  i}^{(\lambda)})\otimes \phi_i^{(\lambda),n}(w_j^{(\lambda)}).  
\nan
It is clear from this that $\zeta$ is determined by the
sequences of linear homomorphisms $\phi^{(\lambda),n}_i$. Note that
\ban ||\zeta_{n+m}-\zeta_{n}|| \to 0, &\quad\quad& n\to\infty
.  
\nan 
Since the $S(t_{ij}^{(\lambda)})$ are linearly independent, 
this implies that for each $\lambda\in\cP_+$ and $i,
j\in\{1, 2, ...,d_\lambda \}$, $\phi^{(\lambda),n}_i(w^{(\lambda)}_j)$
is a Cauchy sequence in $V_\mu$. But since $\Hom (W(\lambda), V_\mu)$
is a finite dimensional complex vector space with the basis $\{
v_r\otimes {\bar w}_j^{(\lambda)}\}$, it is clear that this further
implies that, for each $\lambda\in\cP_+$ and $i\in\{1, 2,
...,d_\lambda \}$, $\phi^{(\lambda),n}_i$ is a Cauchy sequence in
$\Hom (W(\lambda), V_\mu)$ and so
\ban 
\lim_{n\to\infty} \phi^{(\lambda),n}_i = 
\phi^{(\lambda)}_i\in {\rm Hom}_{\C}(W(\lambda), V_\mu).
\nan 

Now we will now show that this limit $\phi_i^{(\lambda)}$ must in fact
be a $\Up$-module homomorphism. 
The defining
property of $\co_q(V_\mu)$ states that 
\ban 
p\circ\zeta =( id_{\ca}\otimes S(p)) \zeta, \quad \forall p\in\Up.
\nan 
Thus, for each $p$, 
$$
||p\circ\zeta_n -( id_{\ca}\otimes S(p)) \zeta_n||\to 0.
$$
Again using the linear independence of the $S(t^{(\lambda)}_{ij})$'s, 
we see that this
implies that, for each $i,k\in \{1,\cdots ,d_\lambda \}$, 
$$
\sum_{j} t_{j k }^{(\lambda)}(S(p)) \phi_i^{(\lambda),n}(w_j^{(\lambda)})
-S(p) \phi_i^{(\lambda),n}(w_k^{(\lambda)})
$$ 
is a null sequence. Thus in the limit we have
\ban
 \phi_i^{(\lambda)}(p w_j^\lambda)
&=& p \phi_i^{(\lambda)}(w_j^\lambda), \quad  \forall p\in\Up.           
\nan  
This is precisely the statement that
the $\phi_i^{(\lambda)}$ are $\Up$-module homomorphisms, 
\ban 
\phi_i^{(\lambda)}&\in&\mbox{Hom}_{\Up}\left(W(\lambda), V_\mu\right)  
\subset \mbox{Hom}_{\C}\left(W(\lambda), V_\mu\right),  
\quad \forall i\in 1,\cdots ,d_\lambda. 
\nan  
It immediately follows from (\ref{Hom}) that  
\ban  \phi^{(\lambda)}_i &=& c_i \, \phi^{(\lambda)},\quad 
c_i\in \C, \nan  
and $\phi^{(\lambda)}$ may be nonzero only when
\ban
{\bar \lambda}&=&\tilde\mu. 
\nan
Hence, if $-\tilde\mu\not\in\cP_+$, we have $\co_q(V_\mu)=0$. 
When $-\tilde\mu\in\cP_+$,  we set $$\nu=(-\tilde\mu)^\dagger.$$ 
Then, we may conclude that $\co_q(V_\mu)$ is spanned by
\ba
\zeta_i &=& \sum_{j}
 S(t_{j i}^{(\nu)})\otimes \phi^{(\nu)} (w_j^{(\nu)}),   \label{irrep}
\na
which are  obviously linearly independent. Furthermore,
\ban x\cdot\zeta_i &=&
 \sum_{j} t^{(\nu)}_{j i} (x)\  \zeta_j,
\ \ \ \ \ x\in\Uq.    \nan
Thus $\co_q (V_\mu)\cong W(\nu)$, 
and this completes the proof of the theorem.  
\end{proof}

We wish to point out that a possible quantum Borel-Weil theorem 
for quantum $GL(n)$ in an algebraic setting (without topology)
without the framework of quantum homogeneous vector bundles 
was elucidated to by Parshall and Wang \cite{Parshall} and 
Noumi et al \cite{Noumi}. Also in \cite{Zhang}, 
a quantum Borel-Weil theorem for the covariant and 
contravariant tensor representations of quantum $GL(m|n)$ 
was obtained along a similar line as that adopted here  
but in an algebraic setting.  
We should also mention that coherent states of compact quantum groups
were investigated in \cite{Jurco} from a representation
theoretical viewpoint. The results reported in that reference
acquire a natural interpretation within the framework of
quantum homogeneous vector bundles.

There are several immediate corollaries of the quantum 
Borel-Weil Theorem \ref{BW}, which are of considerable 
interest. First we note that the proof explicitly constructed the 
isomorphism of the theorem, i.e., (\ref{irrep}). 
 
\noindent {\bf Corollary 1}: {\em
If $\nu^\dagger = -\tilde{\mu}\in\cP_+$, then the following
composition of maps defines the $\Uq$ module isomorphism
$W(\nu)\cong \co_q( V_\mu)$,  }
\ba
W(\nu)
\stackrel{(S\otimes id)P\delta}{\longrightarrow}
\co_q( W(\nu))
\stackrel{id\otimes \phi^{(\nu)}}{\longrightarrow}
\co_q( V_\mu),
\na
{\em where $\phi^{(\nu)}$ is the projection $W(\nu)\rightarrow V_\mu$.}\\
 
Recall that in classical geometry, any analytic function on $G_{\C}/P$
is constant, as the homogeneous space is a  compact complex manifold.
A similar result holds in the quantum  homogeneous space setting.\\
 
\noindent{\bf Corollary 2}:
\ban \begin{array}{l l l}
\Oq(\C)& =& \C \epsilon.
\end{array}\nan
{\em Proof}: This immediately follows from the $\mu=0$ case of
the theorem.

Combining the Corollaries with Proposition \ref{projective},
we obtain\\
 
\noindent{\bf Corollary 3}: {Let $W$ be any finite dimensional $\Uq$-module}.
Then, as $\Uq$-modules,
\ban
\co_q(W) &\cong & \epsilon\otimes W.
\nan

\section{\normalsize APPENDIX:  THE CLASSICAL CASE} \label{classical}

\subsection{Some generalities on homogeneous structures} \label{gen}

Let $G$ be a (real or complex) Lie group and $H$ any
subgroup. Corresponding to each representation $\rho$ of $H$ on a
vector space $V$ one obtains a homogeneous bundle ${\cal V}=
(G\times_H V\to G/H )$ the total space of which is $G\times_H V$, that is
$G\times V$ factored by the equivalence relation
$$
(g,v)=(gh,\rho(h^{-1})v).
$$
Sections of $\cv$ are functions
$$
f:G\to V
$$
satisfying the homogeneity condition
$$
f(gh)=\rho(h^{-1})f(g) .
$$
We will use the notation $\ce\cv$, or $\ce(\cv)$, to mean the sheaf of
germs of smooth sections of $\cv$. By a slight abuse of notation we
will also use this notation to mean simply local sections. 

The space of global smooth sections $\Gamma\ce\cv$ is a $G$-representation
under the action of {\em left translation}, given by
$$
f\mapsto g\cdot f ~~ {\rm for}~~g\in G ~~{\rm and}~~f\in\Gamma\ce\cv ,
$$
where
$
g\cdot f\in \Gamma\ce\cv
$
is the left translated section,
$$
g\cdot f(g')=f(g^{-1}g') ~~\mbox{for all}~~g'\in G.
$$

If $W$ carries a representation $\mu$ of $G$ then the homogeneous bundle
$$
\cw:=G\times_H W
$$
is trivial. The mapping giving
$$
(G/H)\times W\cong G\times_H W
$$
is
$$
(gH,\tilde{w})\leftrightarrow (g,w)~~ {\rm where}~~\tilde{w}=\mu(g)w .
$$
It is easily checked that this is well defined. 
It follows that local sections can be identified
\begin{equation}\label{isom}
\ce((G/H)\times W)\cong\ce(G\times_H W)
\end{equation}
by
$$
\ce((G/H)\times W) \ni  \tilde{f} \leftrightarrow f \in \ce(G\times_H W) 
$$
where $\tilde{f}(g)=\rho(g)f(g)$.
In particular $G\times_H W$ has preferred  sections
of the form $\rho(g^{-1}) w$, where here $w$ is the constant
section of $(G/H)\times W$ corresponding to $w\in W$. Note that
$\Gamma\ce(G\times_H) W$ and $\Gamma\ce((G/H)\times W)$ are each
$G$-representations in two different ways; under the left translation
action of
$G$ and by left multiplication of $\rho(g)$ for $g\in G$. 
It is easily verified that \nn{isom} is {\em not} an isomorphism of
$G$-representations using any combination of these. Note, however,
that if $f(g)=\rho(g^{-1}) w$, for $w\in W$, then 
$$
g'\cdot f(g)=f((g')^{-1}g)=\rho(g^{-1})\rho(g') w
$$ 
and so we have the following proposition.
\begin{proposition}\label{Gisom}
The isomorphism \nn{isom} restricts to give an $G$-monomorphism
$$
W\hookrightarrow \Gamma\ce(G\times_H W)
$$
with where, on the left
hand side, $G$ acts via $\rho$ and, on the right hand side, $G$ acts via
left translation.
\end{proposition}

Recall that as a vector space the Lie algebra $\lig$ of $G$ is just
the tangent space to $G$ at the identity, $\lig=T_e G$. This tangent
space is then identified with the space of left invariant vector
fields via
$$
(Xf)(g):=\left[\frac{d}{dt}f(ge^{tX})\right]_{t=0},
$$
for differentiable functions $f$ on $G$. The space of left
invariant vector fields is closed under commutation and the Lie
bracket on $\lig$ is defined to agree with the commutator of the
elements regarded as left invariant vector fields.  This is consistent
with the adjoint action of $G$ on $\lig$.  Regarding $\lig$ as a
$G$-representation, in this way, allows us to identify the tangent
bundle $TG$ with the homogeneous bundle $G\times_H \lig$.  Note then
that the usual identification of $\lig$ with the right invariant
vector fields is a $G$-monomorphism $\lig\hookrightarrow
\Gamma\ce(G\times_H \lig)$ exactly as in the proposition above.
 
Functions on $G/H$ may clearly be identified with functions on $G$
that are annihilated by $X\in \lih$, where $\lih$ is the Lie algebra
of $H$ identified with the appropriate subalgebra of $\lig$. It
follows then that the tangent bundle $T(G/H)$ to $G/H$ is the homogeneous
bundle,
$$
G\times_H \frac{\lig}{\lih} ,
$$
where the representation of $H$ on $\lig/\lih$ arises from the adjoint
action of $H$ on $\lig$. Note that the right invariant vector fields
on $G$ determine global sections (although not generally globally
non-vanishing) of $T(G/H)$ via the bundle morphism $T(G)\to T(G/H)$ 
determined by the projection $\lig\to \lig/\lih$. 

\subsection{The Setup} \label{setup}

We are interested in studying homogeneous structures on certain
complex homogeneous spaces of  the form $\gc/Q$ where these groups
are described below.
For each complex Lie algebra $\lig$ we write $\lig^{\Bbb R}$
to mean the same Lie algebra but regarded as a
real Lie algebra. 
\newcommand{\strutt}{\rule[-7pt]{0pt}{7pt}}
$$\begin{array}{rcl} 
G^{\Bbb C}&=&\parbox[t]{300pt}{Connected and
simply connected semi-simple complex Lie group. Lie algebra
$\lig$.\strutt}\\ 
Q&=&\parbox[t]{300pt}{A parabolic subgroup of
$G^{\Bbb C}$ with Levi factor $L$.  Lie algebra $\liq$ with
Levi decomposition $\liq=\lil\oplus\liu$.\strutt}\\ 
\theta&=&\parbox[t]{300pt}{A Cartan involution on $\lig^{\Bbb R}$
such that 
$(\theta (\liq^{\br})\cap \liq^{\br})=\lil^{\br}$ and which
fixes a compact real form $\lig_0$ of $\lig$.
Then $\liq$ is $\theta$-stable \cite{KV} and
$\lig=\bliu\oplus\lil\oplus\liu$ where {\em
bar} denotes 
conjugation with respect to $\lig_0$.\strutt}\\ 
G&=&\parbox[t]{300pt}{The real subgroup of $\gc$ with Lie algebra 
$\lig_0$.\strutt}\\
K&=&\parbox[t]{300pt}{$G\cap Q$. Lie algebra $\lik=\lig_0\cap
\lil$. \strutt}\\
\end{array}$$
Note that given $\gc$ and $Q$ one can find a maximal toral subalgebra
of $\lig$ and corresponding root decomposition such that $Q$ is a
standard parabolic. In terms of this root decomposition it is an
elementary exercise to describe $\theta$
explicitly. In particular $\theta$ exists. 
(Note that the setup described above is the classical
analogue of the quantum situation with which the article except here
we have allowed $\gc$ to be semisimple rather than restricting to the
simple case.)

By construction $G$ is a compact real form of $\gc$.
In view of our connectivity assumptions it follows that
$K_{\Bbb C}$ is the Levi part $L$ of $Q$.
Observe that, since $Q$ is closed, $K$ is compact. 
Note also that we clearly have a natural inclusion
$$
G/K\hookrightarrow \gc/Q.
$$
Now dim$(\lig_0/\lik)$=dim${}_{\Bbb R}(\lig/\liq)$ so $G/K$ is an open $G$-orbit in $\gc/Q$. On the other hand
$G/K$ is compact so,  
\begin{equation}\label{ident}
G/K=\gc/Q.
\end{equation}
This identification shows that the symmetric space
$G/K$ is naturally endowed with the structure of a complex manifold.
Our starting point is $\cq:=\gc/Q$ and we are interested in
using this identification to analyse the complex homogeneous bundles
on this in terms of real structures on the left hand side.  
Meanwhile note that it follows immediately from the above
that any element $g\in\gc$ may be written
$$
g=g_0q
$$
where $g_0\in G\subset \gc$ and $q\in Q$. Here $g_0$ is any element of
$G$ such that $g_0K$ corresponds to $gQ$ under the identification of
\nn{ident}. Of course $g_0$ is only determined up to right
multiplication by elements of $K$. To be precise we have 
$$
\gc=G\times_K Q,
$$
where the right hand side means $G\times Q$ modulo the equivalence
relation \mbox{$(g_0,q)$ $\sim$ $(g_0k,k^{-1}q)$} for $k\in K$.
We will describe $g_0q$ as a $GQ$-decomposition of
$g\in\gc$. 
  
We wish to study homogeneous bundles on the homogeneous space $\cq$.
Let $\rho$ be a complex representation of $Q$ on a complex vector
space $V$ and denote by $\cv_Q$ the corresponding induced bundle on
$\gc/Q$. By restriction $\rho$ gives a representation of $K$ on $V$.
Let us also denote this representation by $\rho$ and denote by $\cv$
the homogeneous bundle on $G/K$ induced by this representation.  In
view of \nn{ident} both $\cv_Q$ and $\cv$ are natural structures on
$\cq$.  We will show that, regarded as a real structure, $\cv_Q$ may
be identified with the ($K$-induced) homogeneous bundle $\cv$.

First observe that, clearly,
$$
G\times V\hookrightarrow \gc\times V .
$$
If we compose this set inclusion mapping with the onto mapping to
equivalence classes,
$$
\gc\times V \to \gc\times_Q V, 
$$ 
we obtain a mapping
$$
G\times V \to \gc\times_Q V.
$$
Now since $K=G\cap Q$ this clearly factors through the surjective equivalence
mapping $G\times V \to G\times_K V$. That is there is a natural embedding of
the total space manifolds 
$$
G\times_K V \hookrightarrow\gc\times_Q V .
$$
On the other hand let $(g,v)$ be a representative of any element of $\gc\times_Q
V$. According to our observation above there is an element
$g_0\in G$ and an element $q\in Q$ such that $g=g_0q$. Thus
$(g_0,\rho(q)v)$ is an element of $G\times V\subset \gc\times V$
representing the same equivalence class as $(g,v)$. It follows that 
$$
G\times_K V =\gc\times_Q V ,
$$
as claimed. It follows immediately that we may identify the spaces
of sections
\begin{equation} \label{real}
\Gamma\ce(G\times_K V)\cong \Gamma\ce(\gc\times_Q V).
\end{equation}
It is useful to describe this isomorphism explicitly. First observe
that a section $v\in \Gamma\ce(\gc\times_Q V)$ determines a section
$\tilde{v}\in \Gamma\ce(G\times_K V)$ by restriction to $G\subset
\gc$. On the other hand given $\tilde{v}\in \Gamma\ce(G\times_K V)$ we
can construct the corresponding function on $\gc=G\times_K Q$ by
\begin{equation}\label{extend}
v(g):= \rho(q^{-1})\tilde{v}(g_0) 
\end{equation}
where $g_0q$ is a $GQ$-decomposition of $g\in\gc$. It is easily
verified that this is invariant under the equivalence $(g,q)\sim
(gk,k^{-1}q)$ and satisfies $v(gq')=\rho({q'}^{-1})v(g)$ for $q'\in Q$.

Since $K$ acts reductively on $\lig_0$ this representation splits as a
$K$-module
$$
\lig_0=\lik\oplus\lip .
$$
Thus, as a $K$-representation, $\lig_0/\lik=\lip$ and so the tangent bundle
is induced from the adjoint action of $K$ on $\lip$.
It follows that for any homogeneous
bundle $\cv$, induced by a representation $\rho$ of $K$, there  are 
natural connections $\nd$, 
$$
\nd : \Gamma\ce\cv\to \Gamma\ce(\cv\otimes T^*(G/K)) 
$$
given by
$$
(\nd_x f) (g)=\langle(\nd f)(g) ,x \rangle :=(x f)(g)+\mu(x)f(g),
$$
for $x\in\Gamma\ce(TG)$, and where $\mu$ is any function on $\lig_0$
taking values in Aut$(V)$ such that
$$
\mu |_{\lik}=\rho |_{\lik} ~~{\rm and}~~
\mu(Ad(k^{-1})X)=\rho(k^{-1})\mu(X)\rho(k) ~~{\rm if }~~k\in K
~{\rm and}~~ X\in \lig_0.
$$
(Here, as throughout this appendix, we use the same symbol, $\rho$
in this case, to represent a representation of a group and the
corresponding derivative representation of the Lie algebra.)  
In certain circumstances there is a natural
choice of the linear function $\mu$. For 
example if the representation $\rho$ of $K$ on $V$ is the restriction
of a representation (that we will also denote $\rho$) of $G$ on $V$
then it is natural to take $\mu=\rho$ as then $\nd$ annihilates the
``constant'' sections $\rho(g^{-1})v$, $v\in V$. (See section
\ref{gen} above).

Since $\cq$ has a complex structure it is natural to extend $\nd$ to an
operator
$$
\nd :\Gamma\ce\cv \to \Gamma\ce(\cv\otimes_{\Bbb C} {\Bbb C}T^*(\cq))
$$
by complex linearity. That is $\nabla_x$ is defined as above where now
$x$ is a section of the complexified tangent bundle ${\Bbb
C}T(\cq)$. As mentioned above, we are interested in the case that the
homogeneous bundle over $\cq$ arises from a $Q$-representation $\rho$
on a complex vector space $V$. In this case we have a semi-natural
definition of $\nd$ where we take $\mu |_{\liq}=\rho |_{\liq}$. If
$\rho$ extends to a $\gc$-representation on $V$ then we also require
$\mu |_{\bliu}=\rho |_{\bliu}$, otherwise we simply take $\mu
|_{\bliu} =0$.  Henceforth $\nd$ refers to such a connection.

Since $\cq$ is complex, ${\Bbb C}T^*(\cq)$ splits into
$(1,0)$ and $(0,1)$ parts, 
$$
\mbox{${\Bbb C}T^\ast(\cq)=T_{1,0}^\ast(\cq)\oplus T_{0,1}^\ast(\cq)$}.
$$
Thus the connection `splits' correspondingly, that is
$$
\nd=\nd^{1,0}\oplus \nd^{0,1}
$$
where $\nd^{1,0}:\Gamma\ce\cv\to\Gamma(\ce\otimes T_{1,0}^\ast(\cq))$
and $\nd^{0,1}:\Gamma\ce\cv\to\Gamma(\ce\otimes T_{0,1}^\ast(\cq))$.
The curvature of the connection $\nd$ on $\cv$ is the section $R$ of
$(\bigwedge^2 {\Bbb C}T^*(\cq))\otimes (\cv\otimes \cv^\ast)$ defined by
$$
R(x, y)v=[\nd_x,\nd_y]v -\nd_{[x, y]}v 
$$
for $x, y\in {\Bbb C}T(\cq)$ and $v\in \ce\cv$.

Now as a $K$-representation ${\Bbb C}\otimes\lip=\lig/\lil$ 
decomposes 
$$
\frac{\lig}{\lil}=\bliu\oplus\liu
$$
and this splitting corresponds precisely to the decomposition of
complex tangent vectors into $(1,0)$- and $(0,1)$-parts.
Consider the curvature of a $\cv$ in the case that
$V$ does not extend to a $\gc$-representation. 
Now, $[\bliu,\bliu]\subset \bliu$, 
thus one obtains immediately from the
definition of $\nd_x$ that, for any $v\in \ce\cv$, if
 $x, y\in\bliu$ then
$$
(R(x, y)v)(g)=(x y v)(g)- (y x v)(g)-([x, y]v)(g)=0 ~~\forall ~~g\in G
$$
since, recall, $[x, y]$ is precisely the left invariant vector field
which acts as $(x y-y x)$. 
Now if $x, y\in\liu$, then, using that $[\liu,\liu]\subset \liu$ and that
$(x\rho(y)v)(g)=(\rho(y)x v)(g)$, we obtain
\ban 
(R(x, y)v)(g)& =& (x y v)(g)- (y x v)(g)+[\rho(x),\rho(y)]v(g)\\ 
              &-&([x, y]v)(g)-\rho([x, y])v(g)=0, \quad \forall g\in G
\nan  
for the same reason as the previous case and also using that $\rho$ is
a $\liq$ representation.
On the other hand, if $x\in \liu$
and $y\in\bliu$ then, since $\mu(y)=0$ and $[\liu,\bliu]\subset \lil$, 
$$
(R(x, y)v)(g)=(x y v)(g)- (y x v)(g)=([x, y]v)(g)=-\rho([x, y])v . 
$$ 
In particular this shows that the curvature is type $(1,1)$. It follows
that $(\nd^{0,1})^2 v=0$ for any smooth section $v$ of $\cv$. That is 
$\cv$ is a holomorphic vector bundle with 
$$
\dbar=\nd^{0,1}.
$$  

By a similar calculation one easily obtains that $R=0$ for any vector
bundle $\cw$ induced from a $\gc$-representation $W$. Thus again
the induced bundle $\cw$ admits a holomorphic structure and 
$\dbar=\nd^{0,1}$. Note that for $x\in\liq$ and $f\in\Gamma\ce\cw$,
$$
x(\rho(g)f(g))=\rho(g)\dbar_x f(g).
$$
Thus the holomorphic sections of $\cw$ correspond to the holomorphic
$W$-valued functions on $\cq$. Using proposition \ref{Gisom} and that
$\cq$ is compact we have then the following.
\begin{proposition}\label{holisom}
There is an isomorphism of
$G$-modules 
$$
W\cong \Gamma\co\cw
$$
with $G$-action as given in the proposition \ref{Gisom}.
\end{proposition}

According to \nn{extend} there is no restriction on the
$K$-homogeneous section $\tilde{v}$ for it to extend to a
$Q$-homogeneous function on $\gc$. The discussion there is implicitly
treating the underlying manifold of the group $\gc$ as a real structure.
By construction the function $v$, obtained from $\tilde{v}$ as in
\nn{extend}, satisfies the identity
\begin{equation}\label{diffextend}
x v+\rho(x)v=0 ~~{\rm for~all}~~x\in\liq ,
\end{equation}
this following immediately from the derivative of \nn{extend}.
This is no condition on $\tilde{v}$ as each $x\in \liu\subset\liq$
gives a (real) left invariant vector field on $\gc$ which, at
$G\subset\gc$, is not tangent to $G$. 

The complex structure of $\cq$ arises from regarding $\gc$ as a
complex manifold and in this case the story is rather different. As
mentioned above it is natural to complexify the tangent space in this
case whence each left invariant vector field $x\in\liu$ along
$G\subset \gc$ can be written as complex linear combination of 
tangent vectors to $G$. Thus, at $G\subset\gc$, \nn{diffextend} becomes an 
equation on $\tilde{v}$,
$$
x\tilde{v}+\rho(x)\tilde{v}=0 ~~{\rm for~all}~~x\in\liq .
$$ 
In other words a
section of $G\times_K V$, that is a function
$$
\tilde{v}: G\to V ~~{\rm such~that}~~
\tilde{v}(g_0k)=\rho(k^{-1})\tilde{v}(g_0)~~{\rm for}~~g_0\in
G,~k\in K ,
$$ 
can only extend to a function
$$
v: \gc\to V ~~{\rm such~that}~~v(gq)=\rho(q^{-1})v(g)~~{\rm for}~~g\in
\gc,~q\in Q 
$$
if 
$$
\dbar v=0. 
$$
 On the other hand if $\tilde{v}$
satisfies \nn{diffextend} on $G$ then we can integrate to obtain a
function of the form $\rho(q^{-1})\tilde{v}(g_0)$, which, as discussed above,
is a section of $G\times_Q V$. Thus $\dbar \tilde{v}=0$ is also a sufficient
condition for $K$-homogeneous functions on $G$ to extend to
$Q$-homogeneous functions on $\gc$.  In summary then:
\begin{proposition} \label{autohol}If we regard
$\cv$ as a holomorphic bundle on the complex manifold $\cq$ then
smooth sections of $\cv$ extend to $Q$-homogeneous functions on $\gc$
if and only if they are holomorphic.
\end{proposition}

\subsection{The Borel-Weil theorem, Projectivity and Frobenius Reciprocity} \label{BWandP}

We are now poised to give an elementary treatment of the Borel-Weil
theorem.
\begin{theorem}
Let $V$ be an irreducible finite dimensional $Q$-representation.
If $\cv$ denotes the homogeneous bundle on $\cq$ induced by $V$. Then,
as a $\gc$-representation,
$$
\Gamma\co\cv=\left\{\begin{array} {cl}
W & \mbox{ if there is a $Q$-epimorphism } W\to V \\
0 & \mbox{ otherwise, }
\end{array} \right.
$$
where $W$ is an irreducible finite dimensional $\gc$-representation.
\end{theorem}
\begin{proof}{Proof}
First note that if $\Gamma\co\cv\neq 0$ then it contains a section $f$
such that $f(e)\neq 0$ (or else by left translation
$\Gamma\co\cv=0$). Thus evaluation at the identity determines a
non-trivial $Q$-homomorphism $\Gamma\co\cv\to V$. Since $V$ is
$Q$-irreducible this is a
surjection. Now as $\cq$ is compact it follows from elliptic theory
that $\Gamma\co\cv$ is finite dimensional. Since $\gc$ is reductive
this representation decomposes and there is an irreducible component
$W$ in $\Gamma\co\cv$ such that there is a $Q$-epimorphism $W\to V$.
 
Next we observe that given such a $Q$-epimorphism $\pi: W\to V$ then,
as $\gc$-representations, $W\hookrightarrow \Gamma\ce\cv$.  Let
$\tilde{\pi}$ be the induced $\gc$-homomorphism
$$
\tilde{\pi}:\Gamma\co\cw \to \Gamma\co\cv
$$
given by
$$
(\tilde{\pi} \tilde{f})(g):=\pi(\tilde{f}(g))
$$
for $g\in \gc$ and $\tilde{f}\in\Gamma\co\cw$. 
Note that it follows from proposition \ref{autohol} that $(\tilde{\pi}\tilde{f})$ is holomorphic since  $\pi$ is a $Q$-homomorphism.
Now recall proposition \ref{holisom} there is a
$\gc$-isomorphism 
$$
W\cong \Gamma\co\cw
$$
and each section $\tilde{f}\in\Gamma\co\cw$ is of the form
$\tilde{f}(g)=\rho(g^{-1})w$ for some $w\in W$. It follows that for
any such non-vanishing section $\tilde{f}$ and given any $w'\in W$ there is
some $g\in\gc$ such that $\tilde{f}(g)=w'$. It follows immediately that
$\tilde{\pi}$ is injective and so $W$ is a $\gc$-submodule of
$\Gamma\co\cv$ as claimed.

Finally we show that $\tilde{\pi}$ is
onto.  Since $\Gamma\co\cv$ is finite dimensional it follows, by the
usual theory of weights, that any element of $\Gamma\co\cv$ is in the
$\gc$-orbit of an element that is annihilated by the left action of
$\bliu$. Since $W$ is irreducible, it suffices to show that a such
$\ol{U}$-invariant section $f$ is in the image of
$\tilde{\pi}$. Now $\ol{U}$ acts on $\gc/Q$ and the orbit of the base point
is an open dense set in $\gc/Q$. By continuity it follows that $f$ is
determined by its value $f(e)$ at the identity $e\in\gc$. 
Now there is an element $w\in W$ such that $\pi(w)=f(e)$
and $\rho(X)w=0$ for all $x\in\bliu$. Let
$\tilde{f}(g):=\rho(g^{-1})w$. Then $\tilde{\pi}\tilde{f}$ is $\ol{U}$
invariant and $(\tilde{\pi}\tilde{f})(e)=f(e)$. Thus
$f=\tilde{\pi}\tilde{f}\in\tilde{\pi}W$ as required to be
shown. 

In summary we have that, if $\Gamma\co\cv\neq 0$, then
there is an irreducible $\gc$-representation $W$ such that there is a
$Q$-epimorphism $W\to V$, and 
$$
\Gamma\co\cv\cong W
$$
as $\gc$-representations. This proves the theorem.
\end{proof}

{\bf Projectivity.} The Borel-Weil theorem, as above, is usually
stated in terms of weights since the condition that there be a
non-trivial $Q$ (or equivalently $\liq$) epimorphism $W\to V$ ($W$ and
$V$ as in the theorem) is equivalent to $V$ having a highest weight
which is dominant for $\lig$. The quantum version of this theorem,
presented in section \ref{BWsec} (theorem \ref{BW}), is expressed in
this manner. Under a weaker condition on the highest weight of $V$ we
can establish, in a natural way, the projectivity of the $\ce$-module
$\Gamma\ce\cv$.

Suppose now, then, that there is a $K$-module epimorphism
$$
\phi:W\to V.
$$ 
Let $V^{\perp}:={\rm ker}\phi$. Then since $K$ is reductive it follows
that as a $K$-module $W=V\oplus V^\perp$. Let $\cv^\perp$ be the
homogeneous bundle over $\cq$ induced by the $K$-representation
$V^\perp$. Regarding $\cw$ and $\cv$ also as $K$-induced bundles, it
is at once clear that $\cw$ may be expressed as a bundle direct sum
$\cw=\cv\oplus \cv^\perp$. It is easily seen that this carries over to
the $\ce$-module of local smooth sections and the $\Gamma\ce$-module of
global smooth sections:
$$
\ce\cw=\ce\cv\oplus \ce\cv^\perp~~{\rm
and}~~\Gamma\ce\cw=\Gamma\ce\cv\oplus \Gamma\ce\cv^\perp .
$$
Now recall that $\cw$ is a trivial bundle with fibre $W$ and so, in
particular, we have 
$$
\Gamma\ce\cv\oplus \Gamma\ce\cv^\perp=\Gamma\ce\otimes_{\Bbb C} W .
$$
This establishes the following theorem.
\begin{theorem} \label{swannlike} Suppose that $V$ is an irreducible $K$-module such
that there is a $K$-module epimorphism, $W\to V$, where
$W$ is an irreducible $G$-module. Then
$\Gamma\ce\cv$ is a finite type projective module over the algebra
$\Gamma\ce$ of functions on $\cq$.
\end{theorem}
Of course this theorem is just a special case of Swann's
theorem. However with a view to establishing the corresponding result
in the quantum case it is useful to expose, as we have, the mechanics
underlying the result.  Note also that we can regard $V$ and $W$ as
modules for the complexified groups and their Lie algebras. In this
picture the condition on $V$ is equivalent to the existence of an
$\lil$-epimorphism $W\to V$. This in turn is equivalent to requiring
that the highest weight for $V$ be in the Weyl group orbit of a
$\lig$-dominant weight. That is that this highest weight be integral
for with respect to $\lig$.

{\bf Frobenius Reciprocity.} An early observation in the proof of the
Borel-Weil theorem above was that if $\Gamma\co\cv\neq 0$ then there
is a $L$-epimorphism $\Gamma\co\cv \to \cv$. The same argument shows
that there is always a  $L$-epimorphism  ${\rm Ev}:\Gamma\ce\cv \to
\cv$ determined by evaluation at the identity. It is an elementary exercise 
to verify that if there is a $\gc$-monomorphism $\iota: W \to
\Gamma\ce\cv$ then the composition ${\rm Ev}\circ \iota$ is a
$L$-epimorphism $W\to V$. Here $\gc$ and its subgroup $L$ act on
$\Gamma\ce\cv $ by left translation.

On the other hand if we have $\phi\in {\rm Hom}_L (W,V)$ then this
determines a $\gc$-homomorphism (with respect to the left translation
action of $\gc$) \newline
\mbox{$\Gamma\ce(\gc\times_L W)  \to\Gamma\ce(\gc\times_L V)$}
in the obvious way. Composing this with the $\gc$-equivariant injection 
$W\hookrightarrow \Gamma\ce (\gc\times_L W)$ described in proposition 
\ref{Gisom} we obtain a $ \gc$-monomorphism 
$\tilde{\phi}: W\to \Gamma\ce(\gc\times_L V)$. 

It is easily verified that ${\rm Ev}\circ \tilde{\phi}=\phi$ and that  
conversely, $\widetilde{Ev\circ \iota}=\iota$. This is just the usual
result of Frobenius reciprocity.
\begin{theorem}\label{classFR}
Let $W$ and $V$ be respectively  $\gc$ and $L$ modules and
let $\cv$ denote the homogeneous bundle induced by $V$. Then
there is a canonical isomorphism 
$$
{\rm Hom}_{\scriptsize G^{\Bbb C}} (W,\Gamma\ce\cv)\cong {\rm Hom}_{\scriptsize G^{\Bbb C}} (W,V),
$$
where $\gc$-action on $\Gamma\ce\cv$ is by left translation.  
 \end{theorem}

\subsection{Geometry, Analysis and Algebra}

Let $G$ be a compact group.
Write $K$ to mean either the field ${\Bbb R}$, of real numbers,
or the field ${\Bbb C}$, of complex numbers.  
There is a left action of $G$ on $K$-valued functions given by {\em right
translation} (c.f.\ left translation discussed above)
$$
f\mapsto g\circ f~~~{\rm where}~~~g\circ f (h)=f(hg)~~g,h\in G.
$$
A $K$-valued
{\em representative function} $f$ is a continuous function
$$
f: G\to K
$$
such that the span of the $G$-orbit of $f$, under right translation,
is finite dimensional. That is the representative functions are just
the functions which generate the finite dimensional $G$-invariant
subspaces of the continuous functions $C(G,K)$ on $G$. We write 
$\ct (G,K)$ to denote the space of representative functions.
(See, for example, \cite{BD} for an introduction to these functions
and their role in the theory Lie groups and their representations.)

Suppose that $\rho$ is a finite-dimensional  representation of $G$
then the matrix elements satisfy
$$
\rho_{ik}(hg)=\sum_j \rho_{ij}(h)\rho_{jk}(g) ,
$$
demonstrating that the right translate of the function $\rho_{ik}$ is a linear
combination of the finite set of functions $\{ \rho_{ij}\}$ arising from
the matrix elements of the representation $\rho$.
Thus the matrix elements of finite-dimensional continuous
representations over $K$ are examples of representative functions. It
is well known that, conversely, when all irreducible representations
are used, such matrix element functions
generate $\ct (G,K)$ as a $K$-vector space. Now considering the matrix
elements of dual, direct sum and tensor product of representations quickly
reveals that the representative functions $\ct (G,K)$   admit a natural $K$-algebra
structure which, as a subalgebra of $C(G,K)$, is closed under complex
conjugation. 

>From our point of view the importance of these special functions
arises from the next two theorems. The proof of these theorems (also
see \cite{BD}) involves some standard analysis and functional
analysis.
 \begin{quotation}The first theorem, which is the celebrated theorem of Peter
and Weyl, states that any continuous or $L^2$ function on $G$ can be
approximated by representative functions:
Recall that a compact Lie group $G$ admits a unique normalised (left
and right) invariant Haar integral. Functions are said to be $L^2$ if
they are square-integrable with respect to this integral.
\begin{theorem}(Peter-Weyl)
The representative functions are dense in both $C(G,{\Bbb C})$ and 
$L^2(G,{\C})$.
\end{theorem}

The second theorem, which is
a cornerstone of Tannaka-Kre\u{\i}n duality theory, implies that the
space of representative functions contains all the information of the
compact Lie group:
Let $G_{\Bbb R}$ be the set of $\br$-algebra homomorphisms $\ct
(G,\br)\to \br$. Each $g\in G$ determines an {\em evaluation
homomorphism} $e_g :\ct(G,\br)\to \br$ by $t\mapsto t(g)$.
\begin{theorem}
The duality map 
$$
i:G\to G_{\br} ,~~\mbox{\em given by}~~ g\mapsto e_g,
$$  
is an isomorphism of Lie groups.
\end{theorem}
\end{quotation}
Thus the analysis behind these theorems provides the faithful link
between the algebra of functions on the group and the group itself,
which we regard as a fundamentally geometric object. 

These results suggest a programme where one studies the Lie group, its
associated structures and representations via the algebra of
representative functions or their completion to say $L^2$-functions.
For example, the classical results discussed in sections \ref{gen} to
\ref{BWandP} above could all be reworked in this picture. (See also
\cite{IS} where they show that for a closed subgroup $H$ of a compact
$G$, $G/H$ may be identified with the algebra of homomorphisms $\ct
(G/H,{\Bbb C})\to \br$, where $\ct(G/H,{\Bbb C})\subset\ct(G,{\Bbb C})$
is the subring of functions which factor through the map $G\to G/H$.)
In the case of quantum groups we are without an underlying concrete
group manifold, but these classical results suggest an approach to
defining and studying analogous geometric notions in terms of
representative type functions on the quantized universal enveloping
algebra.


\small
\begin{thebibliography}{9999}
\bibitem{Manin} Yu. I. Manin, {\em Quantum Groups and Noncommutative
         Geometry}, Universite de Montreal, Centre de Recherches
         Mathematiques, Montreal, PQ (1988). 
\bibitem{Woronowicz} S. L. Woronowicz, 
       {\em Differential calculus on compact matrix
       pseudo groups (quantum groups)},
       Commun. Math. Phys. {\bf 122} (1989) 125. 
\bibitem{Zumino} 
C. S. Chu, P. M. Ho and B. Zumino, {\em Some complex quantum manifolds
       and their geometry}, Preprint (1996).
\bibitem{Majid} S. Majid, {\em Advances in quantum and braided geometry}, 
         Preprint (1996).  
\bibitem{groups} S. L. Woronowicz, {\em Compact matrix pseudo  groups},
       Commun. Math. Phys. {\bf 111} (1987) 613. 
\bibitem{Brzezinski} T. Brzezinski and S. Majid, {\em Quantum group 
      gauge theory on quantum spaces}, Commun. Math. Phys. {\bf 157}
      (1993) 591. Erratum. {\bf 167} (1995) 235.
\bibitem{Durdevic} M. Durdevic, {\em Geometry of principle bundles I}, 
        Commun. Math. Phys. {\bf 175} (1996) 457.  
\bibitem{Connes} A. Connes, {\em Noncommutative geometry}, Academic 
      Press (1994).  
\bibitem{Beastwood} R. J. Baston and  M. G. Eastwood,
               {\em The Penrose Transform; its interaction with
               representation theory}, Oxford University Press, Oxford, (1989).
\bibitem{Jimbo} M. Jimbo, {\em A $q$-difference analogue of 
        of $U({\frak g})$ and the Yang-Baxter equation}, 
        Lett. Math. Phys. {\bf 10} (1985) 63.  
\bibitem{Chari} V. Chari and A. Pressley, {\em A guide to 
        quantum groups}, Cambridge University Press, Cambridge (1994).    
\bibitem{Gould} Y. Z. Zhang and M. D. Gould, {\em Unitarity and complete 
         reducibility of certain modules over quantized affine Lie 
         algebras},  J. Math. Phys. {\bf 34} (1993) 6045. 
\bibitem{Faddeev} L. D. Faddeev, N. Yu. Reshetikhin and L. A. Takhtajan,
       {\em Quantization of Lie groups and Lie algebras},  
       Leningrad Math.  J. , {\bf 1} (1990) 193.
\bibitem{Ruan} E. G. Effros and Z. Ruan, {\em Discrete quantum groups I. 
               The Haar measure}, Internat. J. Math., in press.  
\bibitem{Koorwinder} M. S. Dijkhuizen and T. H. Koorwinder,
       {\em CQG algebras: a direct algebraic approach to compact 
       quantum groups}, Lett. Math. Phys. {\em 32} (1994) 315. 
\bibitem{Schneider} H. J. Schneider,
        {\em Principle homogeneous spaces for arbitrary Hopf algebras},
        Israel J. Math. {\bf 72} (1990) 196.
\bibitem{Lakshmibai} V. Lakshmibai and N. Yu. Reshetikhin, 
       {\em Quantum deformation of flag and Schubert schemes}, 
      C. R. Acad. Sci. Paris. Ser. I. Math. {\bf 313} No 3. 
      (1991) 121-126. 
\bibitem{Dijkhuizen} M. S. Dijkhuizen and T. H. Koorwinder, 
         {\em Quantum homogeneous spaces, duality and quantum 2-spheres}, 
         Geom. Dedicata {\bf 52} (1994) 291. 
\bibitem{Montgomery} S. Montgomery, {\em Hopf algebras and their 
        actions on rings}, Regional Conference Series in Math.,     
        {\bf 82} (1993).
\bibitem{Parshall} B. Parshall and J. P. Wang,
         {\em Quantum linear groups},
         Memoirs Amer. Math. Soc., {\bf 89} No. 439 (1991) 1 - 157.
\bibitem{Noumi} M. Noumi, H. Yamada and K. Mimachi, {\em Finite-dimensional
         representations of the quantum group $GL_q(n,\, \C)$ and
         zonal spherical functions on $U_q(n-1)\backslash U_q(n)$},
         Japanese J. Math., {\bf 19} (1993) 31.
\bibitem{Zhang}  R. B. Zhang, {\em Structure and representation of the
         quantum general supergroup}, Univ. of Adelaide Preprint (1996).
\bibitem{Jurco} B. Jurco and P. Stovicek, {\em Coherent states for 
          compact quantum groups}, Commun. Math. Phys. {\bf 182} 
         (1996) 221.  
\bibitem{BD} T. Bro\"{o}cker and T.\ tom Dieck, {\em Representations
of Compact Lie Groups}, Springer-Verlag, New York (1985).
\bibitem{IS} N.\ Iwahori and M.\ Sugiura,
    {\em A duality theorem for homogeneous manifolds of compact Lie groups},
    Osaka J.\ Math.\ {\bf 3} (1966) 139--153.
\bibitem{KV} A.W.\ Knapp and D.A.\ Vogan, {\em Cohomology, Induction
     and Unitary Representations}, Princeton University Press, Princeton
     (1995).
\end{thebibliography}
\end{document}